\documentclass[a4paper]{article}

\usepackage{chemformula} 
\usepackage[T1]{fontenc} 
\usepackage{amssymb}
\usepackage{amsthm}
\usepackage{amsmath}
\usepackage{mathtools}
\usepackage{subcaption}
\usepackage{color}
\usepackage{tikz}
\usepackage{pgfplots}
\usepackage{upgreek}
\usetikzlibrary{calc,patterns}
\usetikzlibrary{positioning} 
\usetikzlibrary{calc}
\usepackage{pgfplots}
\pgfplotsset{width=10cm,compat=1.9}
\usepackage{multirow}
\usepackage{rotating}
\usepackage{authblk}
\usepackage{siunitx}
\pgfplotsset{every tick label/.append style={font=\large}}

\usepackage{hyperref}
\usepackage[all]{hypcap}
\usepackage[a4paper, total={6.5in, 8.5in}]{geometry}

\newcommand{\leb}{\mathrm{leb}}
\newcommand{\bn}{\boldsymbol{n}}
\newcommand{\bx}{\boldsymbol{x}}
\newcommand{\bi}{\boldsymbol{i}}
\newcommand{\bk}{\boldsymbol{k}}
\newcommand{\bs}{\boldsymbol{s}}
\newcommand{\bA}{\mathrm{\textbf{A}}}
\newcommand{\bB}{\mathrm{\textbf{B}}}
\newcommand{\bC}{\mathrm{\textbf{C}}}
\newcommand{\bL}{\mathrm{\textbf{L}}}
\newcommand{\bG}{\mathrm{\textbf{G}}}
\newcommand{\bF}{\mathrm{\textbf{F}}}
\newcommand{\bQ}{\mathrm{\textbf{Q}}}
\newcommand{\bP}{\mathrm{\textbf{P}}}
\newcommand{\bS}{\mathrm{\textbf{S}}}
\newcommand{\omegaij}{\omega_{ij}^{\eta}}
\newcommand{\ylm}{Y_{\ell m}}
\newcommand{\ylmdash}{Y_{\ell ' m'}}
\newcommand{\clomo}{C_{0\ell_0m_0}^j}
\newcommand{\cik}{C_{\ti\tk}^j}
\newcommand{\fni}{f_n^i}
\newcommand{\ti}{\mathtt{i}}
\newcommand{\tk}{\mathtt{k}}
\newcommand{\xni}{\boldsymbol{x}_i^n}
\newcommand{\xnj}{\boldsymbol{x}_j^n}

\title{Linear scaling computation of forces for the domain-decomposition linear Poisson--Boltzmann method}
\author{Abhinav Jha \footnote{Universit\"at Stuttgart, Institute of Applied Analysis and Numerical Simulation, Pffafenwaldring 57, 70569, Stuttgart, Germany, Email: \texttt{abhinav.jha@ians.uni-stuttgart.de}},
Michele Nottoli\footnote{Universit\"at Stuttgart, Institute of Applied Analysis and Numerical Simulation, Pffafenwaldring 57, 70569, Stuttgart, Germany, Email: \texttt{michele.nottoli@ians.uni-stuttgart.de}},
Aleksandr Mikhalev\footnote{RWTH Aachen University, Applied and Computational Mathematics, Schinkelstra\ss e 2, 52062, Aachen, Germany, Email: \texttt{mikhalev@acom.rwth-aachen.de}},
Chaoyu Quan\footnote{SUSTech International Center for Mathematics, and Guangdong Provincial Key Laboratory of Computational Science and Material Design,  Southern University of Science and Technology, Shenzhen, China,
Email: \texttt{quanchaoyu@gmail.com}},
Benjamin Stamm\footnote{Universit\"at Stuttgart, Institute of Applied Analysis and Numerical Simulation, Pffafenwaldring 57, 70569, Stuttgart, Germany, Email: \texttt{benjamin.stamm@ians.uni-stuttgart.de}}}
\date{}

\begin{document}
\maketitle

\begin{abstract}
The Linearized Poisson--Boltzmann (LPB) equation is a popular and widely accepted model for accounting solvent effects in computational (bio-) chemistry. 
In the present article we derive the analytical forces of the domain-decomposition-based ddLPB-method with vdW or SAS surface.
We present an efficient strategy to compute the forces and its implementation,  allowing linear scaling of the method with respect to the number of atoms using the fast multipole method (FMM). 
Numerical tests illustrates the accuracy of the computation of the analytical forces and compares efficiency with other available methods.
\end{abstract}

\section{Introduction}\label{sec:intro}
Most chemical processes and virtually all biochemical processes happen in condensed phase, a situation where the reacting part, or in general the studied part, is embedded in an environment which usually consists of a solvent. For this reason, solvation models, which take into account the effect of the environment on the interesting part (solute), are widely used in computational chemistry and biochemistry.
These models can be broadly divided into two classes, explicit solvation models and implicit (continuum) solvation models.
Explicit solvation models consider the molecular representation of both, the solute and solvent, making the method more accurate, but computationally expensive and also dependent on a large set of empirical parameters (force field).
On the other hand, continuum solvation models treat the solvent as a continuum, described only by a few macroscopic properties. This approach, by its nature, cannot describe specific interactions and anisotropic environment, however it presents some large advantages, it reduces the computational cost significantly, requires fewer parameters and implicitly takes into account the sampling over the degrees of freedom of the solvent. 
For this reason, implicit solvation models are nowadays popular computational approaches to characterize solvent effects in the simulation of properties and processes of molecular systems in condensed phase \cite{TP94,HN95,RS99,CT99,OL00,TMC05}.

Independently from the choice between explicit or implcit solvation model, the solute can be modelled by different levels of theory ranging from (possibly polarizable) force-fields up to coupled cluster theory within a multi-scale approach. 
Thus, this wide scope of different models of the solute has made implicit solvation models popular in different application areas as, depending on the level of theory, structures ranging from only a few atoms to thousands or millions are considered.

The Poisson--Boltzmann (PB) equation is one of the widely used implicit solvation model that we will consider in this paper. The PB equation were described independently by Gouy already in 1910 \cite{gouy_sur_1910} and Chapman in 1913 \cite{chapman_li_1913} to equate the chemical potential and relative forces acting on a small adjacent volumes in an ionic solution between two plates having different voltages. Debye and H\"uckel generalised this concept in 1923 \cite{DH23} by applying it to the theory of ionic solutions leading to a successful interpretation of thermodynamic data. The solutions to the nonlinearised equation were sought by Gronwall, \cite{GL28} in function terms with powers of the inverse of the dielectric constant as coefficients. Simpler electrostatic models for globular proteins were put forward quite early, \cite{Ki34, Li24, NT67}, while DNA and other linear polyelectrolytes were later specialised with proper structural parameter (see \cite{LK54, ABM51, Kat71, Man78}). All the aforementioned models were based around the PB equation or its linear approximation and led to quite accurate results.

We consider here specifically the linearized Poisson--Boltzmann (LPB) equation which describes the electrostatic potential, $\psi$ of the solvation model in the following form
\begin{equation}\label{eq:lpb_equation}
-\nabla\cdot \left[ \varepsilon(\bx)\nabla \psi(\bx)\right]+\overline{\kappa}(\bx)^2\psi(\bx)=4\pi \rho_{\mathrm{M}}(\bx)\qquad \mathrm{in}\ \ \mathbb{R}^3,
\end{equation}
where $\varepsilon(\bx)$ is the space-dependent dielectric permittivity function, $\overline{\kappa}(\bx)$ is the modified Debye-H\"uckel parameter, and $\rho_{\mathrm{M}}(\bx)$ is the solute charge distribution.

\begin{figure}
    \centering
    \includegraphics[width=8cm]{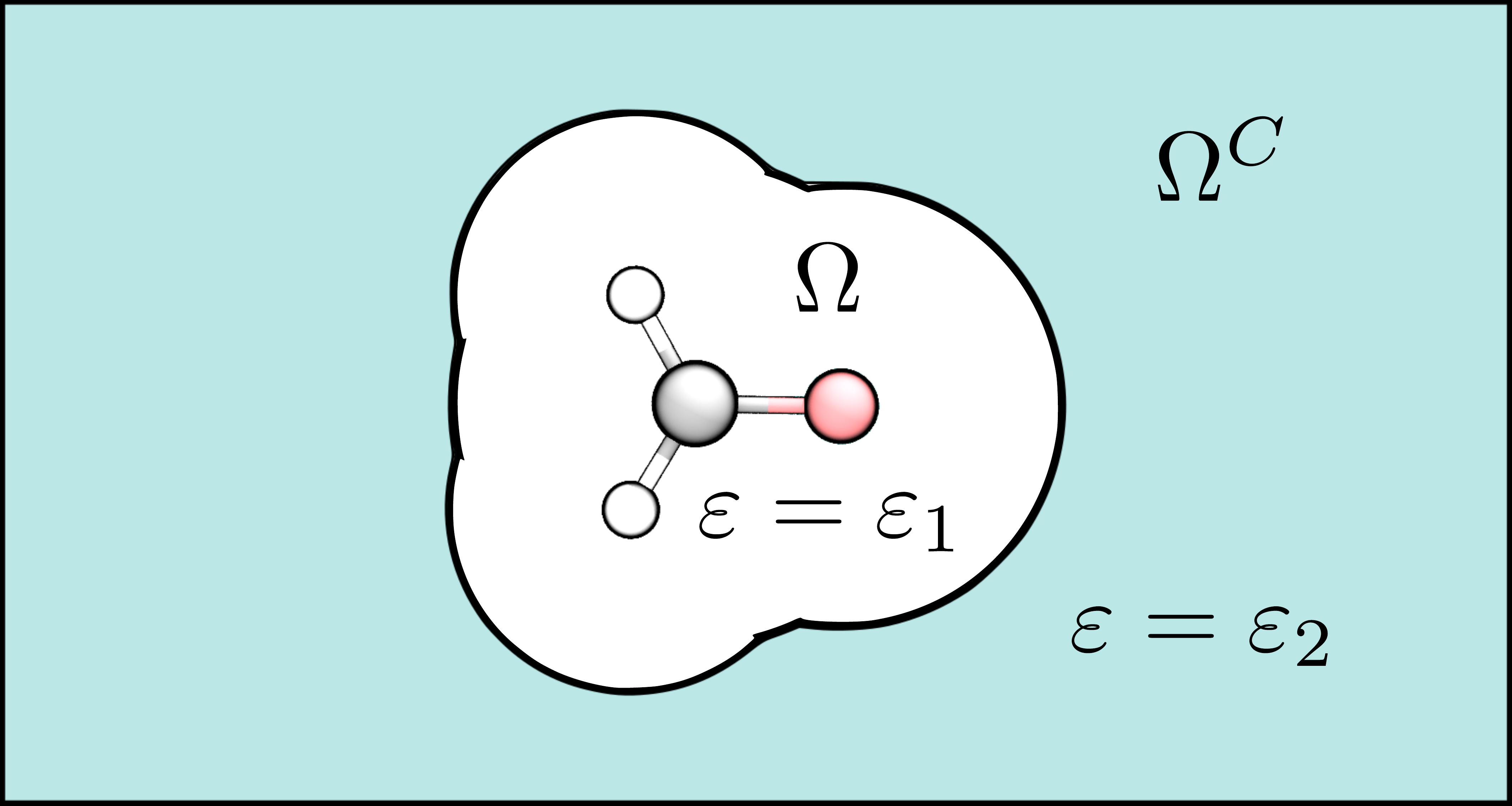}
    \caption{Cavity example for a formaldehyde molecule}
    \label{fig:solute_cavity}
\end{figure}

We denote the solute cavity by $\Omega$ and the solvent region by $\Omega^{\mathrm{C}}=\mathbb{R}^3\setminus \Omega$. To describe the solute-solvent region we will use the van-der Waals (vdW) surface (see Fig.~\ref{fig:solute_cavity}). The solute cavity $\Omega$ is defined as a union of overlapping subdomains, $\Omega_j$, i.e.,
$$
\Omega=\bigcup_{j=1}^M \Omega_j,
\qquad 
\Omega_j=B_{r_j}(\bx_j),
$$
where each $\Omega_j$ is a vdW ball with radius $r_j$ and center $\bx_j$, and $M$ is the total number of atoms. Then $\varepsilon(\bx)$ has the form
\begin{equation*}
\varepsilon(\bx) = \begin{cases} 
        \varepsilon_1 & \mbox{ in } \Omega,\\
        \varepsilon_2  & \mbox{ in } \Omega^{\mathrm{C}},
      \end{cases}
\end{equation*}
where $\varepsilon_1$ and $\varepsilon_2$ are the solute and solvent's dielectric permittivity, respectively. Furthermore, $\overline{\kappa}(\bx)$ has the form
\begin{equation*}
\overline{\kappa}(\bx) = \begin{cases} 
        0 & \mbox{ in } \Omega,\\
        \sqrt{\varepsilon_2}\kappa  & \mbox{ in } \Omega^{\mathrm{C}},
      \end{cases}
\end{equation*}
where $\kappa>0$ is the Debye-H\"uckel screening constant of the solvent.

We would like to mention some of the widely used methods for solving the LPB equation such as the boundary element method (BEM), the finite difference method (FDM), and the finite element method (FEM), and we refer to \cite{LZHM08} for a review. The main idea of the BEM is to recast the LPB equation as an integral equation defined on a two-dimensional solute-solvent interface \cite{YL90, BFZ02, ABWT09, BSR11, SCV22, Cooper2019}. It is an efficient way to solve the LPB equation, which can be optimized using fast multipole methods \cite{ZPH15} and the hierarchial treecode technique \cite{LZHM08}. The PAFMPB solver \cite{LCH10, ZPH15} uses the former optimization technique, whereas the TABI-PB \cite{WK13,Wilson2022} uses the latter one. 
The PB-SAM solver developed by Head-Gordon et al. \cite{lotan2006analytical,yap2010new,yap2013calculating} discretizes the solute-solvent interface (such as the vdW surface) with grid points on atomic spheres like a collocation method and solves the associated linear system by use of the fast multipole method. It primarily targets the interaction of disjoint molecular compounds.
However, one of the limitations of all these solvers relying on integral equations and layer potentials is that it cannot be generalized to solve the nonlinear PB (NPB) equation as opposed to PDE-based methods such as the FDM or FEM.

The finite difference approach is the most popular method to solve linear or nonlinear PB equations. The main idea is to cover the region of interest with a big-box grid and choose different kinds of boundary conditions. Some of the popular software packages using the FDM include UHBD \cite{Madura95}, Delphi \cite{LiLi12}, MIBPB \cite{CCC10}, and APBS \cite{BSJ01, DCLN07, JE17}. One of the drawbacks of the FDM is that the cost can increase considerably with respect to the grid dimension.

The finite element approach, compared to FDM, provides more flexible mesh refinement and a proper convergence analysis \cite{CHX07}. The SDPBS and SMPBS offer fast and efficient approximations of the size-modified PB equation \cite{Xie14, YX15, JXYXY15, XYX17}.

Alongside the PB model in the quantum mechanical (QM) community, continuum solvation models such as the conductor-like screening model (COSMO), proposed in \cite{KS93}; the surface and simulation of volume polarization for electrostatics (SS(V)PE) \cite{Ch99, Ch06}; polarizable continuum model (PCM) \cite{TMC05, MST81, CMT97, BC98, CRSB03} have been developed as a cheap but in a physically sound manner to include solvation effects in the QM description of a molecule and it's properties.  The classical PCM and the COSMO model can be considered as the special cases for PB solvation models. In the classical PCM,  the solvent is represented as a polarizable continuous medium that is non-ionic ($\kappa =0$) whereas the COSMO is a reduced version of the PCM, where the solvent is represented as a conductor-like continuum. Some of the common ways of solving the COSMO model is the BEM \cite{CT95} or the York-Karplus method \cite{YK99}. For the PCM model numerical methods include \cite{CSRB02, SBKPSF04, SF10, LH10a, LH10b}.

In this paper, we focus on the domain decomposition (dd) framework. Recently, in \cite{QSM19} a domain decomposition algorithm has been proposed for the LPB equation, which uses a particular Schwarz domain decomposition method to solve Eq.~\eqref{eq:lpb_equation}. The ideas of the ddLPB method can be traced back to the domain decomposition methods proposed for the COSMO model (ddCOSMO) \cite{CMS13, LSCMM13, LLSS14, LSLS14} and the PCM model (ddPCM) \cite{SCLM16, GLS17, NSSL19}. These methods do not require any mesh or grid of the molecular surface, are easy to implement, and about two orders of magnitude faster than the state of the art \cite{LLSS14}. In particular, the ddCOSMO solver can perform up to thousands of times faster than equivalent existing algorithms.

Similar to the aforementioned dd approaches the ddLPB method does also not require any mesh or grid but depends, as ddCOSMO and ddPCM, only on the Lebedev quadrature points \cite{LL99} on a two-dimensional sphere. Hence it is convenient to be applied in molecular dynamics without re-meshing the molecular surface as is required for the BEM. 
The ddLPB solver adopts a spectral Galerkin method for discretization and benefits from high sparsity of the involved matrices for the Laplace and screening Poisson equations in $\Omega$, which are coupled by a non-local integral equation on the boundary.
The latter takes the majority of cost but can be further accelerated using for example the fast multipole method (FMM).
Numerical implementations show that the ddLPB solver is very efficient even without acceleration techniques (see \cite{QSM19} for details).

The focus of this work is to develop the framework of the computation of first derivatives of the solvation energy with respect to some parameters and the forces in particular,
for the LPB method in the domain decomposition paradigm. 
The electrostatic solvation force is given by the negative gradient of the solvation energy with respect to the nuclear positions and encompasses the reaction field force (RFF), the dielectric boundary force (DBF), and the ionic boundary force (IBF). Out of the three forces, the computation of DBF is quite challenging. The development for computing the DBF was initiated around 30 years ago by Davis and MacCammon in \cite{DM91} where they developed the algorithm based on the Maxwell stress tensor for the two dielectric model which was investigated further through a variational approach in \cite{CDLM08}. At the same time as Davis and MacCammon an alternative algorithm was developed for BEM using the induced surface charge in \cite{Zau91}. Similar results to \cite{Zau91} were obtained using a Maxwell stress tensor for the FEM in \cite{CF97}. Approaches for computing the DBF using FDM were investigated in \cite{GDLM93, IBR98} using a sufficiently smooth-varying dielectric permittivity constant at the molecular surface, but many models used in practise assume a sharp interface. To circumvent this problem a new formulation was proposed on the concept of boundary polarization charge in \cite{CYWL11}. This idea was further expanded to include the abrupt transitional dielectric in \cite{CYL12}. In this work we present the results regarding the total electrostatic solvation forces, which combines RFF, DBF, and IBF 
altogether. 
As can be deduced from above, the computation of the different force components seems to be well-established for FDM and FEM while it seems much harder to generalize this concept for methods based on sharp interfaces, such as the BEM. 
For example, up to our knowledge, the computation of forces is not implemented in popular software such as the TABI-PB method.

Our approach is different and based on the analytical gradients of the discrete energy using the adjoint-method, see, e.g.~\cite{cossi2002,LSCMM13}.
Thus, upon the controllable residual of solving the adjoint linear system, the computed derivatives are the exact negative derivatives of the solvation energy with respect to the nuclear coordinates.

The choice of the solute-solvent interface is part of the model and can be described by the vdW-surface, solvent accessible surface (SAS) or the solvent excluded surface (SES).
For a given solute molecule, both, the SAS and SES, were first introduced by Lee \& Richards in the 1970s~\cite{lee1971interpretation,Richards77} and reflect some properties of the solvent by reducing the solvent molecules to spherical probes~\cite{TMC05}. The SAS is, as the vdW-surface, the surface of a union of balls, but with increased radii compared to the vdW-cavity.
The SES is also called ``the smooth molecular surface'' or ``the Connolly surface'', due to Connolly’s fundamental work~\cite{connolly1983analytical}, and has been rigorously defined and analyzed in~\cite{quan2016mathematical} yielding an efficient meshing-algorithm~\cite{quan2017meshing}. 

It is known that SAS- or vdW-type surfaces yield less accurate energies compared to SES-based models, in particular for large molecular structures. 
On the other hand, SES-based models are not very efficient if high-accurate numerical approximations are required, in particular if high accuracy is needed for the computation of accurate forces or if the PB-model is coupled to quantum-mechanical Hamiltonians. 
We therefore think, as a first step towards a SES-based domain-decomposition methods, that a linear scaling method for the computation of energy and forces arising from the linearized Poisson-Boltzmann equation based on SAS- or vdW-surfaces is a valuable contribution to the state of the art. 

By the nature of the adjoint method, the derivation of the terms in the computation of the forces is a very technical task, but a necessity in order to make accessible the method to models requiring the gradient of the solvation energy with respect to the nuclear coordinates such as molecular dynamics or geometry optimization. 
Further, we accelerate the method based on an FMM-implementation which has recently been proposed in \cite{mikhalev_nottoli_stamm_2022} for the ddPCM model yielding a linear scaling method for the computation of the energy and forces.

The paper is divided as follows: Section~\ref{sec:ddlpb_method} introduces the notations and gives a summary of the domain decomposition algorithm for the LPB equation.  In Section~\ref{sec:computation_of_forces} we derive the adjoint method to compute analytical derivatives and the forces. In Section~\ref{sec:numres} we present a comprehensive numerical study, before we conclude in Section~\ref{sec:summary}. Lastly, in Appendix~\ref{appendix_a} we give the details of the FMM terms needed for the ddLPB method.

\section{ddLPB Method}\label{sec:ddlpb_method}
\subsection{Linear Poisson--Boltzmann Equations}
One notes that the LPB equation~\eqref{eq:lpb_equation} can be written as two equations, one defined in the solute cavity $\Omega$, namely the Laplace equation given by
\begin{equation}\label{eq:laplace_equation}
-\Delta \psi_r=0\qquad \mathrm{in}\ \ \Omega,
\end{equation}
which is obtained from transforming the Poisson equation by using the transformation $\psi_r=\psi-\psi_0$ where $\psi_0$ is the potential generated by $\rho_{\mathrm{M}}$ in the vacuum, i.e.,
\begin{equation}\label{eq:potential_vaccum}
-\Delta \psi_0=\frac{4\pi}{\varepsilon_1}\rho_{\mathrm{M}}\qquad \mathrm{in}\ \ \mathbb{R}^3;
\end{equation}
and a homogeneous screened Poisson (HSP) equation defined on the solvent region given by
\begin{equation*}
-\Delta \psi(\bx)+\kappa^2\psi(\bx)=0\qquad \mathrm{in}\ \ \Omega^{\mathrm{C}},
\end{equation*}
Using potential theory arguments one can define the HSP equation inside $\Omega$,
\begin{equation}\label{eq:hsp_equation}
-\Delta \psi_e(\bx)+\kappa^2\psi_e(\bx)=0\qquad \mathrm{in}\ \ \Omega,
\end{equation}
{with two classical jump conditions}
\begin{align*}
[\![\psi]\!] & =  0\quad \mathrm{on}\ \ \Gamma,\nonumber\\
[\![\partial_{\bn}\left(\varepsilon \psi\right)]\!] & =  0\quad \mathrm{on}\ \ \Gamma, 
\end{align*}
along the solute-solvent interface $\Gamma=\partial\Omega$, and where $[\![f]\!]$ denotes the jump of the function $f$, given by $[\![f]\!]=f|_\Omega - f|_{\Omega^{\mathrm{C}}}$, and $\partial_{\bn}\psi$ the normal derivative of $\psi$.
Based on the classical jump condition of $\psi$, a coupling condition between Eq.~\eqref{eq:laplace_equation} and Eq.\eqref{eq:hsp_equation} arises through a function $h$ defined by
\begin{equation}
\label{eq:nonlocal-coupling}
h=\mathcal{S}_{\kappa}\left(\partial_{\bn}\psi_e-\frac{\varepsilon_1}{\varepsilon_2}\partial_{\bn}\left(\psi_0+\psi_r\right)\right)\qquad \mathrm{on}\ \ \Gamma,
\end{equation}
where $\mathcal{S}_{\kappa}:H^{-1/2}(\Gamma)\rightarrow H^{1/2}(\Gamma)$ denotes a single-layer operator on $\Gamma$ and $H^{\pm 1/2}(\Gamma)$ denote the fractional Sobolev spaces \cite{Ada75}.

We call $\psi_r$ and $\psi_e$ the reaction potential and the extended potential, respectively. In this paper, we assume that the solute's charge distribution $\rho_{\mathrm{M}}$ is supported in $\Omega$ and in particular given by the sum of $M$ point charges, i.e.,
\begin{equation}\label{eq:rho_definition}
\rho_{\mathrm{M}}(\bx)=\sum_{i=1}^Mq_i\delta(\bx-\bx_i),
\end{equation}
where $q_i$ denotes the (partial) charge carried on the $i^{\mathrm{th}}$ atom with center $\bx_i$, and $\delta$ is the Dirac delta distribution, {but the framework can easily be generalized to non-classical charges under the usual assumption $\mathrm{supp}\left(\rho_{\mathrm{M}}\right)\subset \Omega$.}

\subsection{Domain Decomposition Algorithm}
The domain decomposition algorithm that we will consider in this paper has been derived in \cite{QSM19}. For brevity, we will not be deriving the whole method, but we will only present the main equations {required for the derivation of analytical forces.}

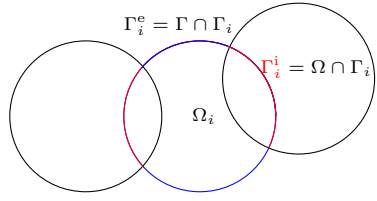
\begin{figure}[t!]
\begin{center}
\begin{tikzpicture}
\draw[blue] (0,0) arc (0:360:1cm);
\draw[black] (-1.5,0) arc (0:360:1cm);
\draw[black] (1.3,0.5) arc (0:360:1cm);
\draw[red] (0,0) arc (0:222:1cm);
\draw[blue] (0,0) arc (0:140:1cm);
\draw[red] (0,0) arc (0:68:1cm);
\draw[red] (0,0) arc (0:-25:1cm);
\begin{scriptsize}
\draw (1.4,0) node[ text width= 5cm, align=left] {$\Omega_i$};
\draw (2.3,0.65) node[ text width= 5cm, align=left] {$\textcolor{red}{\Gamma_i^{\mathrm{i}}}=\Omega\cap \Gamma_i$};
\draw (0.5,1.2) node[ text width= 5cm, align=left] {$\textcolor{black}{\Gamma_i^{\mathrm{e}}}=\Gamma\cap \Gamma_i$};
\draw (-3.4,0) node[ text width= 5cm, align=left] {};
\end{scriptsize}
\end{tikzpicture}
\caption{2-D schematic diagram of $\Gamma_i^{\mathrm{i}}$ and $\Gamma_i^{\mathrm{e}}$}\label{fig:external_internal}
\end{center}
\end{figure}

We first introduce certain notations and functions that will be used throughout the paper. We denote the characteristic function on $\Omega_i$ by $\chi_i$, i.e.,
\begin{equation*}
\chi_i(\bx) := \begin{cases} 
        1 & \mbox{ if } \bx\in \Omega_i,\\
        0  & \mbox{ else},
      \end{cases}
\end{equation*}
and then let
\begin{equation}\label{eq:def_omega_ij}
\omega_{ij}(\bx):=\frac{\chi_j(\bx)}{\sum_{k\in N_i}\chi_k(\bx)}\qquad \mathrm{for}\ \ \bx\in \Gamma_i,
\end{equation}
where $N_i$ denotes the set of indices of spheres intersecting $\Omega_i$ ($i$ not included).  {We make the convention that if $|N_i|=0$, we define $\omega_{ij}(\bx)=0$ for all $j$.} The boundary $\Gamma_i$ of the sphere $\Omega_i$ can either be on the solute-solvent boundary, $\Gamma$, i.e., on the external part or inside the solute cavity, i.e., the internal part. To distinguish between the two cases we define the characteristic function, $\chi_i^{\mathrm{e}}(\bx)$ as
\begin{equation*}
\chi_i^{\mathrm{e}}(\bx) := \begin{cases} 
        1 & \mbox{ if } \bx\in \Gamma_i^{\mathrm{e}},\\
        0  &  \mbox{ if } \bx\in \Gamma_i^{\mathrm{i}},
      \end{cases}
\end{equation*}
where $\Gamma_i^{\mathrm{e}}$ and $\Gamma_i^{\mathrm{i}}$ denote the external and internal part of the boundary $\Gamma_i$ respectively, see Fig.~\ref{fig:external_internal}. 
With the definition of $\omega_{ij}(\bx)$ from Eq.~\eqref{eq:def_omega_ij} we have the relation
\begin{equation}\label{eq:def_chi_ie}
\chi_i^{\mathrm{e}}(\bx)=1-\sum_{j\in N_i}\omega_{ij}(\bx) \qquad \mathrm{for}\ \ \bx\in \Gamma_i.
\end{equation}

We define the radial scaling function of order $\ell$ depending on the $i^{\mathrm{th}}$ atom by
\begin{equation}\label{eq:definition_radial_scaling}
r_{\ell}^i(\bx):=\left(\frac{|\bx-\bx_i|}{r_i}\right)^{\ell}.
\end{equation}
The angular dependency relative to the $i^{\mathrm{th}}$ atom is denoted by
\begin{equation}\label{eq:definition_ylm}
Y_{\ell m}^i(\bx):=Y_{\ell m}\left(\frac{\bx-\bx_i}{|\bx-\bx_i|}\right),
\end{equation}
where $Y_{\ell m}:\mathbb{S}^2\rightarrow \mathbb{R}$ is the real-valued orthonormal spherical harmonic of degree $\ell$ and order $m$.
Moreover, we define the following radial Bessel function by
\begin{equation}\label{eq:definition_of_bessel_func_first}
\bi_{\ell}^i(\bx)\coloneqq\frac{\ti_{\ell}(|\bx-\bx_i|)}{\ti_{\ell}(r_i)},
\end{equation}
where $\ti_{\ell}(\bx)$ is the modified spherical Bessel's function of the first kind. 

Finally, we have integrals over the unit sphere $\mathbb{S}^2$ which will be numerically approximated using the Lebedev quadrature rule \cite{LL99} with $N_{\leb}$ points. The approximation over the sphere $\Omega_i$ is given by
$$
\left\langle f,g\right\rangle_{n,i}:=\sum_{n=1}^{N_{\leb}}\omega_n f(\xni)g(\xni),
$$
where $\xni=\bx_i+r_i\bs_n$, $\bs_n\in \mathbb{S}^2$,  and $\omega_n$ is the quadrature weight.

The fully discretized domain decomposition algorithm for the LPB equation gives rise to the system of equations given by
$$
\bL X=g,
$$
where
\begin{equation}\label{eq:AXg}
\bL=\begin{bmatrix}
\bA & 0 \\
0 & \bB
\end{bmatrix}+
\begin{bmatrix}
\bC_1 & \bC_2 \\
\bC_1 & \bC_2
\end{bmatrix}, \quad
X=\begin{bmatrix}
X_r \\
X_e
\end{bmatrix},\quad \mathrm{and}\quad
g=\begin{bmatrix}
\bG_0+\bF_0\\
\bF_0
\end{bmatrix}.
\end{equation}
The matrices $\bA,~\bB,~\bC_1$, and $\bC_2$ are of the size $M\left(\ell_{\max}+1\right)^2\times M\left(\ell_{\max}+1\right)^2$ where $\ell_{\max}$ denotes the maximum degree of spherical harmonics. The vectors $\bG_0$ and $\bF_0$ on the right-hand side correspond to $\psi_0$ and $\partial_{\bn}\psi_0$, respectively, and $X_r$ and $X_e$ denote the solution vectors corresponding to the reaction potential and the extended potential, respectively. 
After calculating $X$, we can approximate $\psi_{r}$ and $\psi_{e}$ respectively by a linear combination of spherical harmonics as follows
\begin{equation}
\psi_{r}(\bx) \approx \sum_{\ell=0}^{\ell_{\rm max}}\sum_{m=-\ell}^{\ell} [X_{r}]_{i\ell m} \, r_{\ell}^i(\bx) \, Y_{\ell m}^i(\bx), \quad \bx\in \Omega_i,
\end{equation}
and
\begin{equation}
\psi_{e}(\bx) \approx \sum_{\ell=0}^{\ell_{\rm max}}\sum_{m=-\ell}^{\ell} [X_{e}]_{i\ell m} \, \bi_{\ell}^i(\bx) \, Y_{\ell m}^i(\bx), \quad \bx\in \Omega_i.
\end{equation}
\\
We now show the specific formulas of the matrices.
The $(i\ell m, j\ell'm')^{\mathrm{th}}$ matrix entry for  $\bA$ is given by,
\begin{align}\label{eq:definition_of_A_matrix}
\left[ \bA_{ii}\right]_{\ell \ell'}^{mm'} & :=  \delta_{\ell \ell'}\delta_{mm'},\nonumber \\
\left[ \bA_{ij}\right]_{\ell \ell'}^{mm'} & :=  -\sum_{n=1}^{N_{\leb}}\omega_n\omega_{ij}(\xni)r^j_{\ell'}(\xni)Y_{\ell'm'}^j(\xni)Y_{\ell m}(\bs_n),\quad i\neq j,
\end{align}
and the $(i\ell m, j\ell'm')^{\mathrm{th}}$ matrix entry for $\bB$ is given by,
\begin{align}\label{eq:definition_of_B_matrix}
\left[ \bB_{ii}\right]_{\ell \ell'}^{mm'} & :=  \delta_{\ell \ell'}\delta_{mm'},\nonumber \\
\left[ \bB_{ij}\right]_{\ell \ell'}^{mm'} & :=  -\sum_{n=1}^{N_{\leb}}\omega_n\omega_{ij}(\xni)\bi^j_{\ell'}(\xni)Y_{\ell'm'}^j(\xni)Y_{\ell m}(\bs_n),\quad i\neq j.
\end{align}
We note that both the matrices $\bA$ and $\bB$ are sparse in nature, as blocks are nonzero only for interlocking vdW balls.

Next, we move to the matrices $\bC_1$ and $\bC_2$ where the $(i\ell m, j\ell'm')^{\mathrm{th}}$  entry of $\bC_1$ is given by
\begin{align}\label{eq:definition_of_C1_matrix}
\left[ \bC_1\right]_{i\ell m}^{j\ell'm'} & :=  \frac{\varepsilon_1}{\varepsilon_2}\left( \sum_{n=1}^{N_{\leb}}\omega_n\chi_i^{\mathrm{e}}(\xni)Y_{\ell m}(\bs_n)\left[\bQ\right]_{j\ell'm'}^{in}\frac{\ell'}{r_j}\right),
\end{align}
and for $\bC_2$ by,
\begin{align}\label{eq:definition_of_C2_matrix}
\left[ \bC_2\right]_{i\ell m}^{j\ell'm'} & :=  -\left( \sum_{n=1}^{N_{\leb}}\omega_n \chi_i^{\mathrm{e}}(\xni)Y_{\ell m}(s_n)\left[\bQ\right]_{j\ell'm'}^{in}\frac{\ti'_{\ell'}(r_j)}{\ti_{\ell'}(r_j)}\right),
\end{align}
where matrix $\bQ$ is  a matrix of size $M\left(\ell_{\max}+1\right)^2\times MN_{\leb}$ and the $(j\ell'm',in)^{\mathrm{th}}$ entry is given by
\begin{align}\label{eq:definition_of_Q_matrix}
\left[ \bQ\right]_{j\ell'm'}^{in} & :=  \sum_{\ell_0 m_0}\cik \left[ \bP \chi_j^{\mathrm{e}}\right]_{\ell_0m_0}^{\ell'm'} \bk_{\ell_0}^j(\xni)Y_{\ell_0m_0}^j(\xni),
\end{align}
where 
$\bk_{\ell_0}^{j}(\bx)$ is defined similarly to Eq.~\eqref{eq:definition_of_bessel_func_first} given by
\begin{equation}\label{eq:definition_of_bessel_func_second}
\bk_{\ell_0}^{j}(\bx)=\frac{\tk_{\ell_0}(|\bx-\bx_j|)}{\tk_{\ell_0}(r_j)},
\end{equation}
$\tk_{\ell_0}(\bx)$ is the modified spherical Bessel's function of the second kind,
$$
\cik=\left( \frac{\ti'_{\ell_0}(r_j)}{\ti_{\ell_0}(r_j)}-\frac{\tk'_{\ell_0}(r_j)}{\tk_{\ell_0}(r_j)}\right)^{-1},
$$ 
and the notation $\sum_{\ell m}$ denotes $\sum_{\ell=0}^{\ell_{\max}}\sum_{m=-\ell}^{\ell}$.

The matrix $\bP \chi_j^{\mathrm{e}}$ is of size $\left(\ell_{\max}+1\right)^2\times \left(\ell_{\max}+1\right)^2$ whose $\left( \ell_0m_0, \ell'm'\right)^{\mathrm{th}}$ entry is given by
\begin{align}\label{eq:definition_of_P_matrix}
\left[ \bP \chi_j^{\mathrm{e}}\right]_{\ell_0m_0}^{\ell'm'} &:= \sum_{n=1}^{N_{\leb}}\omega_n\chi_j^{\mathrm{e}}(\xnj)Y_{\ell_0m_0}(\bs_n)Y_{\ell'm'}(\bs_n).
\end{align}
Finally, we have the right-hand side vectors.  The $(i\ell m)^{\mathrm{th}}$ entry of the vector $\bG_0$ is given by
\begin{equation}\label{eq:definition_of_G_0}
\left[ \bG_0\right]_{i\ell m}=-\sum_{n=1}^{N_{\leb}}\omega_n \chi_i^{\mathrm{e}}(\xni)\psi_0(\xni)Y_{\ell m}(\bs_n),
\end{equation}
where 
\begin{equation}\label{eq:defintion_of_psi_0}
\psi_0(\bx)=\sum_{j=1}^M\frac{q_j}{\varepsilon_1|\bx-\bx_j|},
\end{equation}
is the solution of Eq.~\eqref{eq:potential_vaccum} and the $(i\ell m)^{\mathrm{th}}$ entry of $\bF_0$ is given by
\begin{equation}\label{eq:definition_of_F_0}
\left[ \bF_0\right]_{i\ell m}=-\frac{\varepsilon_1}{\varepsilon_2}\left(\sum_{n=1}^{N_{\leb}}\omega_n\chi_i^{\mathrm{e}}(\xni)Y_{\ell m}(\bs_n)\sum_{j=1}^M\left[\bS\right]_{jin} \right),
\end{equation}
where
\begin{equation}\label{eq:defintion_of_S}
\left[\bS\right]_{jin}=\sum_{\ell_0 m_0}\cik \clomo\bk_{\ell_0}^j(\xni)Y_{\ell_0m_0}^j(\xni),
\end{equation}
and
\begin{equation}\label{eq:definition_of_c_0_l_0_m_0}
\clomo=\sum_{n=1}^{N_{\leb}}\omega_n \chi_j^{\mathrm{e}}(\xnj)\partial_{\bn}\psi_0(\xnj)Y_{\ell_0m_0}(\bs_n).
\end{equation}

\section{Computation of Forces}\label{sec:computation_of_forces}
The computation of the electrostatic solvation energy, $E_s$ in \cite{QSM19}, follows the ideas of \cite{FBM02} where the reaction potential was used to compute $E_s$. 
For the computation of forces, we require the whole electrostatic potential and hence we define $E_s$ as
\begin{eqnarray}\label{eq:energy_computation}
E_s = \frac{1}{2}\left\langle \psi_{\mathrm r},\rho_{\mathrm M}\right\rangle
 =  \frac{1}{2}\sum_{j=1}^M\left\langle X, Q\right\rangle_j,
\end{eqnarray}
where $X$ is given in Eq.~\eqref{eq:AXg}, $Q$ has the same size as $X$ with
$$
[Q]_{j\ell m}= 
\begin{dcases} 
q_j\delta_{\ell 0}\delta_{m0}, & \mbox{ if } 1\leq j\leq M,\\
0 & \mbox{ if } M<j\leq 2M.
\end{dcases}
$$
and the inner product $\langle \cdot,\cdot\rangle_j$ is given by
$$
\langle X,Q\rangle_j=\sum_{\ell m}\left[X\right]_{j\ell m}\left[ Q\right]_{j\ell m}.
$$

The force with respect to a parameter $\lambda$, such as the position of $\bx_k$ of the $k^{\mathrm{th}}$ atom, is given by,
\begin{equation*}
F_{{\lambda}} = \nabla^{{\lambda}}\left(E_s\right) = \frac{1}{2}\left( \left\langle \nabla^{{\lambda}} X,Q\right\rangle + \left\langle X,\nabla^{{\lambda}} Q\right\rangle\right)=   \frac{1}{2}\left\langle \nabla^{{\lambda}} X ,Q\right\rangle.
\end{equation*}
The ddLPB system is given by $\bL X=g$. Taking the derivative with respect to ${{\lambda}}$:
\begin{eqnarray*}
\nabla^{{\lambda}} \bL X + \bL\nabla^{{\lambda}} X & = & \nabla^{{\lambda}} g \nonumber \\
\nabla^{{\lambda}}X  & = & \bL^{-1}\left(\nabla^{{\lambda}}g-\nabla^{{\lambda}}\bL X\right).
\end{eqnarray*}
Substituting $\nabla^{{\lambda}} X$ in the force computation
\begin{eqnarray*}
F_{{\lambda}}&=& \frac{1}{2}\left\langle \bL^{-1}\left(\nabla^{{\lambda}}g-\nabla^{{\lambda}}\bL X\right), Q \right\rangle\nonumber \\
& = & \frac{1}{2}\left\langle \left(\nabla^{{\lambda}}g-\nabla^{{\lambda}}\bL X\right), \left(\bL^{-1}\right)^*Q \right\rangle,
\end{eqnarray*}
where $\bL^*$ is the adjoint of the matrix $\bL$ and $\left(\bL^{-1}\right)^*Q$ is the solution of the system
\begin{equation}\label{eq:adjoint_system}
\bL^*X_{\mathrm{adj}}=Q.
\end{equation}

Using the definition of $X_{\mathrm{adj}}$ we get the computation of forces as
\begin{equation}\label{eq:force_lpb}
F_{{\lambda}}=\frac{1}{2}\left\langle \left(\nabla^{{\lambda}}g-\nabla^{{\lambda}}\bL X\right), X_{\mathrm{adj}}\right\rangle.
\end{equation}
We note that in Eq.~\eqref{eq:force_lpb} we require the computation of the adjoint system ({but only once for any number of different parameters $\lambda$}) and the derivatives of the $g$ and $\bL$ matrix. The adjoint matrix of the system is given by
\begin{equation}\label{eq:AstarXg}
\bL^*=\begin{bmatrix}
\bA^{\mathrm{T}} & 0 \\
0 & \bB^{\mathrm{T}}
\end{bmatrix}+
\begin{bmatrix}
\bC_1^{\mathrm{T}} & \bC_1^{\mathrm{T}} \\
\bC_2^{\mathrm{T}} & \bC_2^{\mathrm{T}}
\end{bmatrix},
\end{equation}
where $\bA^{\mathrm{T}}$ stands for the transpose of the matrix $\bA$ and respectively others. 

In the next subsection we would present the analytical derivatives that arise in Eq.~\eqref{eq:force_lpb}.

\subsection{Analytical Derivatives}
{We now restrict ourselves to the case where $\lambda$ denotes the central coordinate $\bx_k$ of the $k^{\mathrm{th}}$ atom.} We note that entries of matrix $\bL$ and vector $g$ have certain functions that are not smooth, namely, $\chi_i(\bx)$, $\chi_i^{\mathrm{e}}(\bx)$, and $\omega_{ij}(\bx)$. To define their differentiable counterparts, we follow the ideas presented in \cite{LSCMM13}. We first introduce a polynomial, $p_{\eta}(t)$ given by
$$
p_{\eta}(t):=\eta^{-5}\left(1-t\right)^3\left(6t^2+\left(15\eta -12\right)t+1-\eta^2-15\eta + 6\right),
$$
where $\eta$ is a smoothness parameter. Then the regularized characteristic function is given by
\begin{equation}\label{eq:regular_char_function}
\chi_{\eta}(t) = \begin{dcases} 
        1 & \mbox{ if } t\leq 1-\eta,\\
       p_{\eta}(t)  & \mbox{ if } 1-\eta<t<1, \\
       0 & \mbox{ if }t\geq 1.
      \end{dcases}
\end{equation}

Using Eq.~\eqref{eq:regular_char_function}, the regularized version of $\omegaij(\bx)$ defined in Eq.~\eqref{eq:def_omega_ij} is given by
\begin{equation}\label{eq:regular_omega_ij}
\omegaij(\bx):=d^i(\bx)\chi_{\eta}\left(r_1^j(\bx)\right),\quad \forall \bx \in \Gamma_i,
\end{equation}
with
\begin{equation}\label{eq:di}
d^i(\bx):=\frac{\min \left\lbrace f^i(\bx),1\right\rbrace}{f^i(\bx)},
\end{equation}
where
\begin{equation}
f^i(\bx):=\sum_{k\in N_i}\chi_{\eta}\left(r_1^k(\bx)\right)
\end{equation}
and $r_1^j$ is defined in Eq.~\eqref{eq:definition_radial_scaling}.
Finally, the differentiable counterpart of $\chi^{\mathrm{e}}_i(\bx)$ is given by
\begin{equation}\label{eq:regular_chi_e}
\chi^{\eta}_i(\bx):=1-\sum_{j\in N_i}\omegaij(\bx),\quad \forall \bx \in \Gamma_i.
\end{equation}

One thing to note is that in the definition of $d^i(\bx)$ we have a minimum which is not a smooth function. 
On close inspection we note that if $f^i(\bx)<1$, then $d^i(\bx)=1$, else $d^i(\bx)=1/f^i(\bx)$.


\subsubsection{Sparse Matrices \textbf{A} and \textbf{B}}
\begin{figure}[t!]
\centering
\begin{subfigure}[t!]{0.5\textwidth}
 \centering
 \begin{tikzpicture}
 \draw (0,0) circle (1cm);
 \draw (1.5,0) circle (1cm);
 \begin{scriptsize}
 \draw (0,0) node[ align=center] {$\Omega_i=\Omega_k$};
 \draw (1.5,0) node[ align=center] {$\Omega_{j}$};
 \end{scriptsize}
 \end{tikzpicture}
 \caption{$j\in N_i$ and $i=k$}\label{subfig:i_in_N_j_i_eq_k}
 \end{subfigure}%
\hfill
\begin{subfigure}[t!]{0.5\textwidth}
\centering
 \begin{tikzpicture}
 \draw (0,0) circle (1cm);
 \draw (1.5,0) circle (1cm);
 \begin{scriptsize}
 \draw (0,0) node[ align=center] {$\Omega_i$};
 \draw (1.7,0) node[ align=center] {$\Omega_{j}=\Omega_k$};
 \end{scriptsize}
 \end{tikzpicture}
 \caption{$j\in N_i$ and $j=k$}\label{subfig:i_in_N_j_j_eq_k}
\end{subfigure}%
\vspace{0.2in}
\hfill 
\begin{subfigure}[t!]{0.5\textwidth}
 \centering
 \begin{tikzpicture}
 \draw (0,0) circle (1cm);
 \draw (1.5,0) circle (1cm);
 \draw[black, dashed] (-1.5,0) circle (1cm);
 \begin{scriptsize}
 \draw (0,0) node[ align=center] {$\Omega_i$};
 \draw (1.5,0) node[ align=center] {$\Omega_{j}$};
 \draw (-1.5,0) node[ align=center] {$\Omega_{k}$};
 \end{scriptsize}
 \end{tikzpicture}
 \caption{$j\in N_i$ and $k\in N_i$ and $ k\neq j$}\label{subfig:i_in_N_j_k_in_N_i_k_neq_j}
\end{subfigure}%
\hfill%
 \begin{subfigure}[t!]{0.5\textwidth}
 \centering
 \begin{tikzpicture}
 \draw (0,0) circle (1cm);
 \draw (1.5,0) circle (1cm);
 \draw[black, dashed] (0.75,-1) circle (1cm);
 \begin{scriptsize}
 \draw (0,0) node[ align=left] {$\Omega_i$};
 \draw (1.5,0) node[ align=left] {$\Omega_j$};
 \draw (0.75,-1) node[ align=left] {$\Omega_k$};
 \end{scriptsize}
 \end{tikzpicture}
 \caption{$j\in N_i$ and $k\in N_i$ and $ k\neq j$} \label{subfig:spheres_j_notin_N_i_k_in_N_i_k_in_N_j}
 \end{subfigure}%
\caption{Example of spheres with non-zero contribution in derivatives for $\bA$ and $\bB$.}\label{fig:spheres_sparse_matrices}
\end{figure}
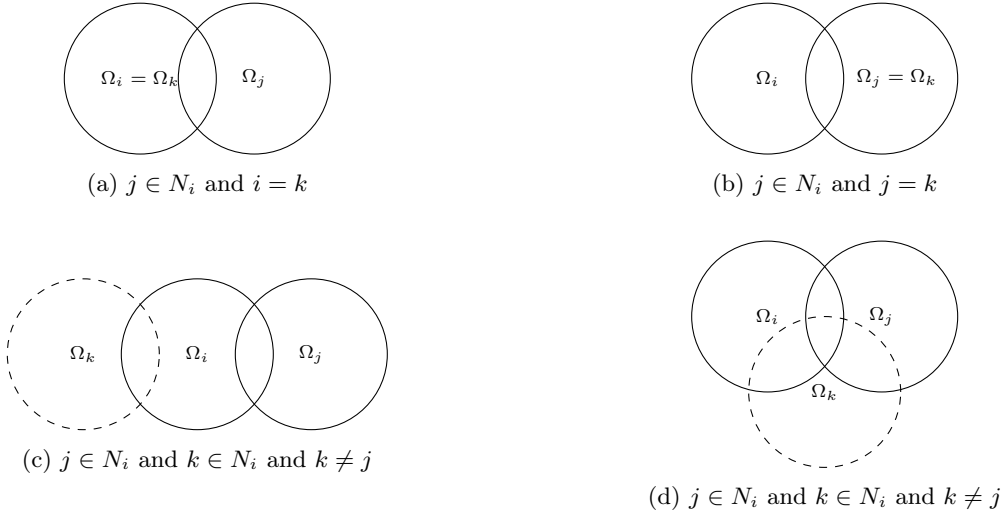

As noted in the previous sections, the matrices $\bA$ and $\bB$ are sparse in nature with constant diagonal entries. As we are finding derivatives with respect to the position of sphere $\Omega_k$, i.e., $\bx_k$, we have the following cases which gives non-zero contribution
\begin{enumerate}
\item $j\in N_i $ and $i=k$ (see Subfig.~\ref{subfig:i_in_N_j_i_eq_k});
\item $j\in N_i $ and $j=k$ (see Subfig.~\ref{subfig:i_in_N_j_j_eq_k});
\item $j\in N_i $ and $k\in N_i$ and $ k\neq j$ (see Subfig.~\ref{subfig:i_in_N_j_k_in_N_i_k_neq_j} and \ref{subfig:spheres_j_notin_N_i_k_in_N_i_k_in_N_j}).
\end{enumerate}

Fig.~\ref{fig:spheres_sparse_matrices} shows the aforementioned cases. Looking at the matrix entries for $\bA$ and $\bB$ we note that we have three terms depending on the position, namely $\omegaij(\xni),\ylmdash^j(\xni),$ and $r_{\ell'}^j(\xni)$ for matrix $\bA$; and $\bi_j^{\ell}(\xni)$ for matrix $\bB$.

For abbreviation, we denote $\nabla_{\bx_k}$ by $\nabla^k$ in the following content.
The derivative of $\omegaij(\xni)$ is given by
\begin{equation}\label{eq:derivative_omegaij}
\nabla^k\left(\omegaij(\xni)\right) = \begin{dcases} 
        d^i(\xni)\left[\chi_{\eta}'\left( r_1^j(\xni)\right)\frac{e^j(\xni)}{r_j}-\delta_{\fni>1}\omegaij(\xni)Z_n^i\right] & \mbox{ if } \left(j\in N_i \right) \wedge \left(k=i\right),\\
        -d^i(\xni)\chi_{\eta}'\left( r_1^j(\xni)\right)\frac{e^j(\xni)}{r_j}
       \left(1-\delta_{\fni>1}\omegaij(\xni) \right)   & \mbox{ if } \left(j\in N_i \right) \wedge \left(k=j\right), \\
       \omegaij(\xni)\delta_{\fni>1}d^i(\xni)\chi_{\eta}'\left(r_1^k(\xni) \right)\frac{e^k(\xni)}{r_k} & \mbox{ if }
        \left(j\in N_i \right) \wedge \left(k\in N_i\wedge k\neq j\right),\\
       0 & \mbox{ else,}
      \end{dcases}
\end{equation}
where
$$
e^j(\bx) \coloneqq ( \bx-\bx_j)/|\bx-\bx_j|,\quad
Z_n^i \coloneqq \sum_{k_0\in N_i}\chi_{\eta}'\left( r_1^{k_0}(\xni)\right)e^{k_0}(\xni)/r_{k_0},
\quad f_n^i \coloneqq f^i(\bx_i^n),
$$
and 
$$
\delta_{\fni>1}=\begin{dcases} 
        1 & \mbox{ if } \fni>1,\\
        0 & \mbox{ if } \fni\leq 1.
        \end{dcases}
$$


Further, the derivative of $\ylmdash^j(\xni)$ is given by 
\begin{equation}\label{eq:derivative_ylm}
\begin{aligned}
    \nabla^k\left(\ylmdash^j(\xni)\right) 
    & =(\nabla\ylmdash)\left(\frac{\bx_i^n-\bx_j}{|\bx_i^n-\bx_j|}\right) \nabla^k \frac{\bx_i^n-\bx_j}{|\bx_i^n-\bx_j|}\\
      & = \begin{dcases} 
        (\nabla\ylmdash)^j\left(\bx_i^n\right)\frac{1}{r_jr_1^j(\xni)} & \mbox{ if }  \left(j\in N_i \right) \wedge (k=i),\\
       -(\nabla\ylmdash)^j\left(\bx_i^n\right)\frac{1}{r_jr_1^j(\xni)} & \mbox{ if }  \left(j\in N_i \right) \wedge (k=j), \\
       0 & \mbox{ else.} 
      \end{dcases}
\end{aligned}
\end{equation}
We now show the details for derivation of Eq.~\eqref{eq:derivative_ylm}. Note that $\forall \bx=(x_1,x_2,x_3),$ we have
\begin{equation}
    \nabla_{\bx} \frac\bx{|\bx|} = \frac1{|\bx|^3}\begin{pmatrix}
|\bx|^2-x_1^2 & -x_1x_2 & -x_1x_3\\
-x_1x_2 & |\bx|^2-x_2^2 & -x_2x_3\\
-x_1x_3 & -x_2x_3 & |\bx|^2-x_3^2
\end{pmatrix}
\end{equation}
and
\begin{equation}
    \nabla_{\bx} \ylmdash\left(\frac\bx{|\bx|}\right)\cdot \bx = 0,
\end{equation}
which yield that 
\begin{equation}
    \nabla_{\bx} \ylmdash\left(\frac\bx{|\bx|}\right)\nabla_{\bx} \frac\bx{|\bx|} = \nabla_{\bx} \ylmdash\left(\frac\bx{|\bx|}\right)\frac1{|\bx|}.
\end{equation}
The equation Eq.~\eqref{eq:derivative_ylm} is then followed. 

Lastly, we have the derivatives of the radial scaling $r^j_{\ell'}(\xni)$ given by
\begin{equation}\label{eq:derivative_radial}
\nabla^k\left(r^j_{\ell'}(\xni)\right) = \begin{dcases} 
        \frac{e^j(\xni)\ell'r_{\ell'-1}^j(\xni)}{r_j}& \mbox{ if } \left(j\in N_i \right) \wedge \left(k=i\right),\\
       -\frac{e^j(\xni)\ell'r_{\ell'-1}^j(\xni)}{r_j}  & \mbox{ if } \left(j\in N_i \right) \wedge \left(k=j\right), \\
       0 & \mbox{ else,}
      \end{dcases}
\end{equation}
and the Bessel scaling, $\bi_{\ell'}^j(\xni)$ which is given by
\begin{equation}\label{eq:derivative_bessel}
\nabla^k\left(\bi_{\ell'}^j(\xni)\right) = \begin{dcases} 
        \frac{\ti'_{\ell'}(|\xni-\bx_j|)e^j(\xni)}{\ti_{\ell'}(r_j)}& \mbox{ if } \left(j\in N_i\right) \wedge \left(k=i\right),\\
       - \frac{\ti'_{\ell'}(|\xni-\bx_j|)e^j(\xni)}{\ti_{\ell'}(r_j)}  & \mbox{ if } \left(j\in N_i\right) \wedge \left(k=j\right), \\
       0 & \mbox{ else.}
      \end{dcases}
\end{equation}

Collecting all the terms, we can compute the derivatives of the $(i\ell m,j\ell'm')^{\mathrm{th}}$ element of matrix~$\bA$ and~$\bB$.
In the case of $i=j$, 
\begin{equation}
\nabla^k\left[ \bA_{ii}\right]_{\ell \ell'}^{mm'} = 0,\quad 
\nabla^k\left[ \bB_{ii}\right]_{\ell \ell'}^{mm'} = 0.
\end{equation}
We now consider different cases for $i\neq j$ as follows.
\begin{enumerate}
\item Case $j\in N_i$ and $k=i$ (Subfig.~\ref{subfig:i_in_N_j_i_eq_k}):
\begin{eqnarray*}
\nabla^k\left[ \bA_{ij}\right]_{\ell \ell'}^{mm'} & = & -\sum_{n=1}^{N_{\leb}}\omega_n\ylm(\bs_n)\Bigg[ r_{\ell'}^j(\xni)\ylmdash^j(\xni)d^i(\xni)\nonumber \\
&& \times \left \lbrace \frac{1}{r_j}\chi_{\eta}'\left( r_1^j(\xni)\right) e^j(\xni)-\delta_{\fni>1}\omegaij(\xni)Z_n^i\right\rbrace\nonumber \\
&& + \omegaij(\xni)\frac{r_{\ell'-1}^j(\xni)}{r_j}\bigg\lbrace \ell'\ylmdash^j(\xni)e^j(\xni) +\left(\nabla\ylmdash\right)^{j}(\xni)\bigg\rbrace\Bigg],
\end{eqnarray*}
and
\begin{eqnarray*}
\nabla^k\left[ \bB_{ij}\right]_{\ell \ell'}^{mm'} & = & -\sum_{n=1}^{N_{\leb}}\omega_n\ylm(\bs_n)\Bigg[ \bi_{\ell'}^j(\xni)\ylmdash^j(\xni)d^i(\xni)\nonumber \\
&& \times \left\lbrace \frac{1}{r_j}\chi_{\eta}'\left( r_i^j(\xni)\right)e^j(\xni)-\delta_{\fni>1}\omegaij(\xni)Z_n^i\right\rbrace\nonumber \\
&& + \omegaij(\xni)\ylmdash^j(\xni)\frac{\ti_{\ell'}'(|\xni-\bx_j|)}{\ti_{\ell'}(r_j)}e^j(\xni) \nonumber \\
&& + \omegaij(\xni)\bi_{\ell'}^j(\xni)\frac{1}{r_jr_1^j(\xni)}\left(\nabla\ylmdash\right)^{j}(\xni)\Bigg].
\end{eqnarray*}
\item Case $j\in N_i$ and $k=j$ (Subfig.~\ref{subfig:i_in_N_j_j_eq_k}):
\begin{eqnarray*}
\nabla^k\left[ \bA_{ij}\right]_{\ell \ell'}^{mm'} & = & \sum_{n=1}^{N_{\leb}}\omega_n\ylm(\bs_n)\Bigg[ r_{\ell'}^j(\xni)\ylmdash^j(\xni)\frac{d^i(\xni)}{r_j}\chi_{\eta}'\left(r_1^j(\xni) \right)e^j(\xni)\nonumber \\
&& \times \left\lbrace 1-\delta_{\fni >1}\omegaij(\xni)\right\rbrace \nonumber \\
&& + \omegaij(\xni)\frac{r_{\ell'-1}^{j}(\xni)}{r_j}\Big\lbrace \ell'\ylmdash^j(\xni)e^j(\xni)+\left(\nabla\ylmdash\right)^{j}(\xni)\Big\rbrace\Bigg],
\end{eqnarray*}
and
\begin{eqnarray*}
\nabla^k\left[ \bB_{ij}\right]_{\ell \ell'}^{mm'} & = & \sum_{n=1}^{N_{\leb}}\omega_n\ylm(\bs_n)\Bigg[ \bi_{\ell'}^j(\xni)\ylmdash^j(\xni)\frac{d^i(\xni)}{r_j}\chi_{\eta}'\left( r_1^j(\xni)\right)e^j(\xni)\nonumber \\
&& \times \left\lbrace 1-\delta_{\fni>1}\omegaij(\xni)\right\rbrace\nonumber \\
&& + \omegaij(\xni)\ylmdash^j(\xni)\frac{\ti_{\ell'}'(|\xni-\bx_j|)}{\ti_{\ell'}(r_j)}e^j(\xni) \nonumber \\
&& + \omegaij(\xni)\bi_{\ell'}^j(\xni)\frac{1}{r_jr_1^j(\xni)}\left(\nabla\ylmdash\right)^{j}(\xni)\Bigg].
\end{eqnarray*}
\item Case  $ \left(j\in N_i \right) \wedge \left(k\in N_i\wedge k\neq j\right)$ (Subfig.~\ref{subfig:i_in_N_j_k_in_N_i_k_neq_j} and~\ref{subfig:spheres_j_notin_N_i_k_in_N_i_k_in_N_j}):
\begin{eqnarray*}
\nabla^k\left[ \bA_{ij}\right]_{\ell \ell'}^{mm'} & = & -\sum_{n=1}^{N_{\leb}}\omega_n\ylm(\bs_n)r_{\ell'}^j(\xni)\ylmdash^j(\xni)\omegaij(\xni)\nonumber \\
&& \times \delta_{\fni>1}d^i(\xni)\chi_{\eta}'\left( r_1^k(\xni)\right)\frac{e^k(\xni)}{r_k},
\end{eqnarray*}
and
\begin{eqnarray*}
\nabla^k\left[ \bB_{ij}\right]_{\ell \ell'}^{mm'} & = &  -\sum_{n=1}^{N_{\leb}}\omega_n\ylm(\bs_n)\bi_{\ell'}^j(\xni)\ylmdash^j(\xni)\omegaij(\xni)\nonumber \\
&& \delta_{\fni>1}d^i(\xni)\chi_{\eta}'\left( r_1^k(\xni)\right)\frac{e^k(\xni)}{r_k}.
\end{eqnarray*}
\end{enumerate}

\subsubsection{Dense Matrices \texorpdfstring{$\bC_1$}{a} and \texorpdfstring{$\bC_2$}{a}}
Now, we move our attention towards the computation of derivatives for {the} matrices $\bC_1$ and $\bC_2$.  We compute the derivative of $\bC_1$ and $\bC_2$ together, i.e., we consider
\begin{equation}
\begin{bmatrix}
\bC_1 & \bC_2\\
\bC_1 & \bC_2
\end{bmatrix}
\begin{bmatrix}
X_r\\
X_e
\end{bmatrix}=
\begin{bmatrix}
\bC_1X_r+\bC_2X_e\\
\bC_1X_r+\bC_2X_e
\end{bmatrix},
\end{equation}
where the $\left(i\ell m\right)^{\mathrm{th}}$ entry of $\left[\bC_1X_r+\bC_2X_e\right]$ is given by:
\begin{eqnarray}\label{eq:def_c1_c2}
\left[\bC_1 X_r + \bC_2 X_e\right]_{i\ell m}& = & \sum_{j=1}^M\sum_{\ell'm'}\sum_{n=1}^{N_{\leb}} \omega_n\ylm(\bs_n)\left[\bQ\right]_{j\ell'm'}^{in}\nonumber \\
&& \times \chi_i^{\eta} \left(\xni\right)\left[ \frac{\varepsilon_1}{\varepsilon_2}\frac{\ell'}{r_j}\left[X_r\right]_{j\ell'm'}-\frac{\ti_{\ell'}'(r_j)}{\ti_{\ell'}(r_j)}\left[X_e\right]_{j\ell'm'}\right]
\end{eqnarray}
We note that we have two terms depending on $\bx_k$, i.e., $\chi_i^{\eta}(\xni)$ and $[\bQ]_{j\ell'm'}^{in}$. Unlike for matrices $\bA$ and $\bB$ we have {non-trivial contributions on the diagonal} as well.  We divide the computation of derivative of Eq.~\eqref{eq:def_c1_c2} into two parts with help of the product rule as follows
\begin{eqnarray*}
\nabla^k \left[\bC_1X_r+\bC_2 X_e\right]_{i\ell m} & = & \sum_{j=1}^M\sum_{\ell'm'}\sum_{n=1}^{N_{\leb}}\omega_n Y_{\ell m}(\bs_n) \left[ \frac{\varepsilon_1}{\varepsilon_2}\frac{\ell'}{r_j}\left[X_r\right]_{j\ell'm'}-\frac{\ti_{\ell'}'(r_j)}{\ti_{\ell'}(r_j)}\left[X_e\right]_{j\ell'm'}\right]\nonumber \\
&& \times \left[ \left[\bQ\right]_{j\ell'm'}^{in}\nabla^k \chi_i^{\eta}(\xni) + \chi_i^{\eta}(\xni) \nabla^k\left[\bQ\right]_{j\ell'm'}^{in} \right].
\end{eqnarray*}

\textbf{Derivative of $\chi_i^{\eta}(\xni)$.}
The first contribution is the derivative of $\chi_i^{\eta}(\xni)$ when keeping $\bQ$ as constant. 
The non zero contribution comes when $k=i$ or $k\in N_i$.  
Combining \eqref{eq:regular_chi_e} and \eqref{eq:derivative_omegaij}, we have
\begin{equation}\label{eq:derivative_chi_i}
\begin{aligned}
    \nabla^k\left(\chi_i^{\eta}(\xni)\right) & = -\sum_{j\in N_i} \nabla^{k}\left(\omegaij(\xni)\right) \\
    & =
    \begin{dcases} 
        \left[(1-\chi_i^{\eta}(\xni))\delta_{\fni>1}-1\right] d^i(\xni) Z_n^i & \mbox{ if } k=i,\\
        \left[1-(1-\chi_i^{\eta}(\xni))\delta_{\fni>1}\right] d^i(\xni)\chi_{\eta}'\left(r_1^k(\xni) \right)\frac{e^k(\xni)}{r_k}
        & \mbox{ if } k\in N_i,\\
        0 & \mbox{ else, }
        \end{dcases}\\
    & =
    \begin{dcases} 
        -\delta_{f_n^i\leq 1} Z_n^i & \mbox{ if } k=i,\\
        \delta_{f_n^i\leq 1}\chi_{\eta}'\left(r_1^k(\xni) \right)\frac{e^k(\xni)}{r_k}
        & \mbox{ if } k\in N_i,\\
        0 & \mbox{ else. }
    \end{dcases}
\end{aligned}
\end{equation}
Here we use the fact that if $f_n^i>1$, then $\chi_i^{\eta}(\xni) =0$; if $f_n^i\leq 1$, then $d^i(\xni) = 1$.

\textbf{Derivative of $\left[ \bQ\right]_{j\ell'm'}^{in}$.}
The second contribution comes from the derivatives of matrix $\bQ$. The entries are given by Eq.~\eqref{eq:definition_of_Q_matrix}. 

In this matrix we note that three terms depend on the position namely, $\left[ \bP \chi_j^{\eta}\right]_{\ell_0m_0}^{\ell'm'}$, $\bk_{\ell_0}^j(\xni)$, and $Y_{\ell_0m_0}^j(\xni)$.  
To be precise, we have
\begin{eqnarray*}
\nabla^k \left[ \bQ\right]_{j\ell'm'}^{in} 
 = &\sum_{\ell_0 m_0}\cik\Big(\nabla^k\left[\bP\chi_j^{\eta}\right]_{\ell_0 m_0}^{\ell'm'}\bk_{\ell_0}^j(\xni)Y_{\ell_0m_0}^j(\xni)\\
& +\left[\bP\chi_j^{\eta}\right]_{\ell_0 m_0}^{\ell'm'} \nabla^k\bk_{\ell_0}^j(\xni)Y_{\ell_0m_0}^j(\xni) + \left[\bP\chi_j^{\eta}\right]_{\ell_0 m_0}^{\ell'm'} \bk_{\ell_0}^j(\xni)\nabla^k Y_{\ell_0m_0}^j(\xni) \Big).
\end{eqnarray*}


The non-zero contribution of the derivative for $\bk_{\ell_0}^j(\xni)$ and $Y_{\ell_0m_0}^j(\xni)$ comes when $k=i$ or $k=j$.
The derivative of $\bk_{\ell_0}^j$ is given by:
\begin{equation}\label{eq:Derivative_of_k_l0}
\nabla^k\left(\bk_{\ell_0}^j(\xni)\right) = \begin{dcases} 
        \frac{\tk'_{\ell_0}(|\xni-\bx_j|)e^j(\xni)}{\tk_{\ell_0}(r_j)}& \mbox{ if } k=i,\\
       - \frac{\tk'_{\ell_0}(|\xni-\bx_j|)e^j(\xni)}{\tk_{\ell_0}(r_j)}& \mbox{ if } k=j,
      \end{dcases}
\end{equation}
while the derivative of $Y_{\ell_0m_0}^j(\xni)$ is already given by Eq.~\eqref{eq:derivative_ylm} with $\ell',m'$ replaced by $\ell_0,m_0$.
The final contribution comes from the derivative of $\left[\bP\chi_j^{\eta}\right]_{\ell_0m_0}^{\ell'm'}$.  We have the computation of
\begin{eqnarray*}
\nabla^k\left[\bP\chi_j^{\eta}\right]_{\ell_0m_0}^{\ell'm'}
=\sum_{{n_0}=1}^{N_{\leb}}\omega_{n_0}\nabla^k\left(\chi_j^{\eta}(\bx^{n_0}_j)\right)Y_{\ell_0m_0}\left(\bs_{n_0}\right)\ylmdash\left(\bs_{n_0}\right),
\end{eqnarray*}
where the derivative of $\chi_j^{\eta}\left(\bx^{n_0}_j\right)$ is given by Eq.~\eqref{eq:derivative_chi_i} with $i,n$ replaced by $j,n_0$.

\subsubsection{Right-hand Side  \texorpdfstring{$\bG_0$ and $\bF_0$}{a}}
The final derivatives we require are those of the right-hand side $\bG_0$ and $\bF_0$. In Eq.~\eqref{eq:definition_of_G_0} we have two terms depending on $\bx_k$; $\chi_i^{\eta}(\xni)$ and $\psi_0(\xni)$. The derivatives of $\chi_i^{\eta}(\xni)$ is given by Eq.~\eqref{eq:derivative_chi_i} and the derivative of $\psi_0$ is given by
\begin{equation}\label{eq:derivative_psi_0}
\nabla^k\left(\psi_0(\xni) \right) = \begin{dcases} 
       \frac{-1}{\varepsilon_1}\sum_{j=1}^Mq_j\frac{e^j(\xni)}{|\xni-\bx_j|^2}& \mbox{ if } k=i,\\
       \frac{q_k}{\varepsilon_1}\frac{e^k(\xni)}{|\xni-\bx_k|^2}  & \mbox{ if } k=j.
      \end{dcases}
\end{equation}
Next we move towards the computation of derivative for $\bF_0$. We note that {the} entries of $\bF_0$ are very similar to the entries of $\left[\bC_1X_r+\bC_2X_e\right]$, with only the addition of the term $\partial_{\bn}\psi_0(\xni)$. The computation of other terms namely, $\chi_i^{\eta}(\xni)$, $\bk_{\ell_0}^j(\xni)$, and $Y_{\ell_0m_0}^j(\xni)$ has been taken before. The derivatives of $\partial_{\bn}\psi_0$ is given by
\begin{equation}\label{eq:derivative_normal_psi_0}
\nabla^k\left(\partial_{\bn}\psi_0(\xni) \right) = \begin{dcases} 
       \sum_{j=1,j\neq i}^M q_j\left[ \frac{3(\xni-\bx_j)(\xni-\bx_j)^{\mathrm{T}}}{|\xni-\bx_j|^5} - \frac{\mathbb{I}_{3\times 3}}{|\xni-\bx_j|^3}\right]\cdot \bn& \mbox{ if } k=i,\\
      - q_k\left[ \frac{3(\xni-\bx_k)(\xni-\bx_k)^{\mathrm{T}}}{|\xni-\bx_k|^5} - \frac{\mathbb{I}_{3\times 3}}{|\xni-\bx_k|^3}\right]\cdot \bn  & \mbox{ if } {k\neq i},
      \end{dcases}
\end{equation}
where $\mathbb{I}_{3\times 3}$ is the identity matrix of size $3\times 3$ and $\bn = \bs_n$ at $\bx_i^n$ is the unit normal derivative.

The computation of forces can be summarized as follows:
\begin{enumerate}
    \item Solve Eq.~\eqref{eq:AXg} to get the reaction potential $X_r$ and the extended potential $X_e$.
    \item Solve Eq.~\eqref{eq:AstarXg} to get the adjoint solution $X_{\mathrm{adj}}$.
    \item Compute the analytical derivatives of the matrix $\bL$ and the right-hand side $g$ with respect to a parameter $\lambda$.
    \item Contract the analytical derivatives with the adjoint solution to get the forces.
\end{enumerate}

\section{Numerical Simulations}\label{sec:numres}
In this section, we present an extensive study for the computation of the electrostatic solvation energy and the electrostatic solvation forces. Before presenting the examples, we would like to mention some details on solving the system of equations~\eqref{eq:AXg} and~\eqref{eq:adjoint_system}.  We follow a slightly different approach as presented in~\cite{QSM19}. We re-write our system of equations~\eqref{eq:AXg} as
$$
\left(\bL_{\bA \bB}+\bL_{\bC}\right)X = g,
$$
where
$$
\bL_{\bA \bB}:=\begin{bmatrix}
\bA & 0 \\
0 & \bB
\end{bmatrix}\quad \mathrm{and}\quad
\bL_{\bC}:=\begin{bmatrix}
\bC_1 & \bC_2 \\
\bC_1 & \bC_2
\end{bmatrix}.
$$

We solve the above system using direct inversion in the iterative subspace (DIIS)~\cite{pulay1980convergence,rohwedder2011analysis} with $\bL_{\bA \bB}$ as the preconditioner.

The initial iterate for this system is taken as $\bL_{\bA \bB}^{-1}g$, i.e. we start with $X=0$.  We refer to these iterations as macro-iterations, and as one needs to solve two linear systems within the preconditionner for finding $\bA^{-1}$ and $\bB^{-1}$, we refer to them as micro-iterations as they are also performed in an iterative manner. The initial iterate for the two linear systems is zero for the first iteration. For the subsequent iterations, we take the solution of the previous macro-iterations as the guess.
Compared to the strategy presented in~\cite{QSM19}, this technique allows our method to be more consistent as one can use the same solver for both the micro and macro-iterations.

For each linear system, the stopping criterion is on the relative increment of the solution, i.e.,
\begin{equation}
\frac{\displaystyle  \|X^{(\nu)} - X^{(\nu - 1)}\|_{\infty}}{\|X^{(\nu)}\|_{\infty}} \leq \mathtt{tol}\qquad \nu\geq 1,
\end{equation}
where $\|\cdot\|_{\infty}$ is the $\ell^{\infty}$-norm of the corresponding vector. However, we use two different tolerances for the micro and macro-iterations, namely, the inner tolerance is equal to the outer tolerance divided by 100.

\begin{table}[ht]
    \centering
    \begin{tabular}{|cccc|}
       \hline
        \textbf{PDB Code} & \textbf{Number of} & \textbf{Name} & \textbf{Reference} \\ 
         & \textbf{Atoms ($\boldsymbol{M}$)}& & \\\hline\hline
        1ay3 & 25 & Nodularin & \cite{A96} \\
        1etn & 180 & Enterotoxin & \cite{O91} \\
        1du9 & 380 & Scorpion toxin & \cite{X00} \\
        1d3w & 2049 & Ferredoxin & \cite{C00} \\
        1jvu & 3964 & Ribonuclease A & \cite{V01} \\
        1qjt & 9046 & EH1 domain & \cite{W99} \\
        1a3n & 10087 & Human haemoglobin & \cite{Tame2000} \\
        1ju2 & 20260 & Hydroxynitrile lyase & \cite{Dreveny2001} \\ \hline
    \end{tabular}
    \caption{Information about the input structures.}
    \label{tab:structures}
\end{table}

The code was tested on a set of input structures with different number of atoms, spanning from $10^1$ to $10^4$ atoms. We prepared the input structures using the tool PDB2PQR provided in the APBS software package\cite{JE17}, the AMBER force field was used to assign the atomic partial charges\cite{Ponder2003}.

The (relative) dielectric constant of the solute's region is set to 1 (vacuum) and the dielectric constant of the environment is set to 78.54 (water). We included two ions of charge $+1$ and $-1$, both in concentration 0.1~M, which combined with a temperature of 298.15~K, correspond to $\kappa = 0.104$.

The radii were assigned in a subsequent step, according to a definition of a solvent accessible surface (SAS): for each atom we set its radius to its value as reported in ref.~\cite{B64} plus a contribution from the effective size of the solvent (1.4~\AA~for water). Table~\ref{tab:structures} reports detailed information about the structures. The same radii were used in the ddX, APBS finite difference method. For what concerns TABI-PB calculation, these can only be done on smooth cavities generated using Nanoshaper. In this case we used the same Van der Waals radii defined in ref~\cite{B64}, but then we generated a solvent excluded surface (SES) using a probe radius of 1.4~{\AA}.

All the calculations were run on the BwUniCluster2.0 using the ``thin'' nodes. These servers are equipped with two Intel Xeon Gold 6230 (2.1~GHz) CPUs, for a total of 80 cores and up to 192~GB of RAM, which run Red Hat Enterprise Linux 8.4 (Ootpa). Furthermore the ddX, APBS and TABI-PB executables were compiled using the Intel compiler 2021.4.0 and linked against the Intel MKL libraries bundled in the same package.

For all the simulations we used 10 cores, except in Sec.~\ref{sec:comparison_software}, where we used a single core while comparing different methods.

For the ddLPB calculations we used our implementation of ddX, available on GitHub \cite{ddX} at commit \texttt{6bbea05} and compiled using the flags \texttt{-O3 -xHost -fp-model=precise}. For the APBS-FDM calculations we used the APBS code available on GitHub (\texttt{Electrostatics/apbs}) at commit \texttt{e8d1a9c} compiled using the default release flags. Finally, for the TABI-PB calculations we used the TABI-PB code available on GitHub (\texttt{Treecodes/TABI-PB}) at commit \texttt{0710ff7} and compiled using the default release flags. TABI-PB also requires the NanoShaper executable, for which the version 0.7.8 was used.

In the following, we will present numerical results, that
 are divided into two parts. We first present the results regarding the accuracy of the method and then we present the results regarding the complexity of the method.

\subsection{Accuracy of the Discretization}

\subsubsection{Numerical Validation of the Analytical Forces}

The analytical forces computed by Eq. \eqref{eq:force_lpb} have been tested against numerical forces that were computed through finite differences. Indeed, the numerical forces are evaluated using the following definition
\begin{equation}
\label{eq:numerical}
   D_{h}[E_s](\lambda) := \frac{E_s (\lambda + h) - E_s (\lambda)}{h} \text{.}
\end{equation}
Here, $\lambda$ is a generic parameter, for instance one component of a nuclear coordinate, and $0 < h \ll 1$ is a small step size.
Note that the ddLPB-method proposed in this manuscript computes the analytical forces, i.e. the exact derivative of the discrete energy, up to the tolerance of the resolution of the adjoint linear system, and the numerical forces are just computed for purpose of testing the former one.

For the numerical test, we selected the two smallest structure (1ay3, 1etn) since the computation of the numerical forces acting on each nuclear coordinate is quite expensive and we computed all the numerical derivatives with respect to the nuclear coordinates using Eq.~\eqref{eq:numerical}, for various finite step sizes. Due to high computational cost related to the repeated number of calculations, we used a coarser discretization: $\ell_{\max}=2$, $N_{\leb}=110$, and $\mathtt{tol}=10^{-6}$. Also, given the small size of the structures we decided to not use the FMM acceleration.

Due to the finite difference approximation of the analytical derivative we expect a first-order convergence of
\[
    \mathsf{Err}_{j,\alpha}(h) := D_{h}[E_s](x_{j,\alpha}) - \frac{\partial E_s}{\partial x_{j,\alpha}},
    \qquad \mbox{with }\;
    \bx_j = ( x_{j,1}, x_{j,2}, x_{j,3})^T,
\]
with respect to $h$. 
As comparison, note that the force acting on the component $\alpha=1,2,3$ of nuclei $j$ due to the solvation model is given by $F_{j,\alpha}= - \frac{\partial E_s}{\partial x_{j,\alpha}}$.

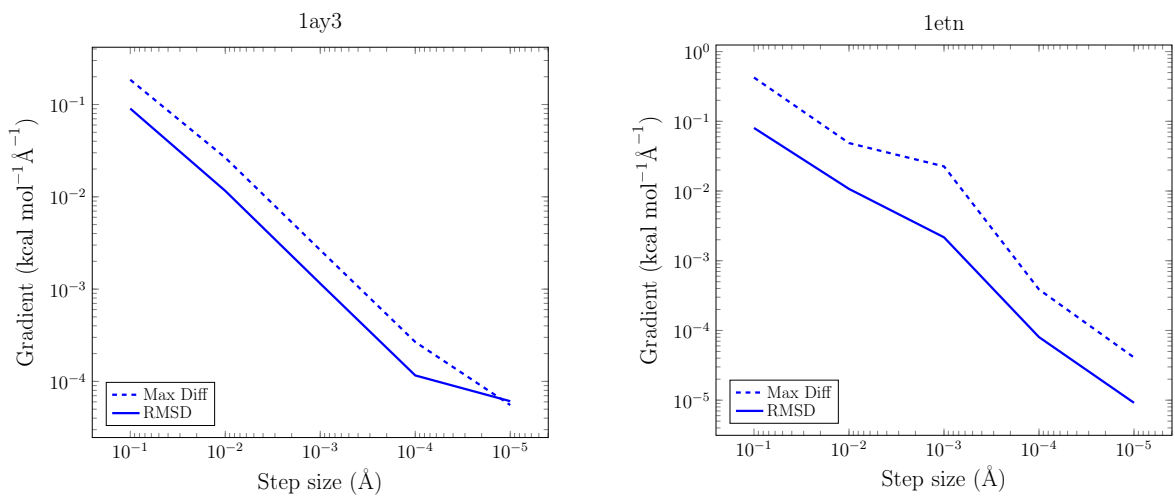
\begin{figure}
\centering
\begin{subfigure}[t!]{0.49\textwidth}
 \centering
 \begin{tikzpicture}[scale=0.6]
 \pgfplotsset{
      scale only axis,
  }
\begin{loglogaxis}[
    legend pos=south west, xlabel = \Large{Step size (\AA)},ylabel= \Large{Gradient ($\mathrm{kcal\ mol}^{-1}\mbox{\AA}^{-1}$)},
    legend cell align ={left},
    x dir=reverse, title = \Large{1ay3}, 
    legend style={nodes={scale=0.75, transform shape}}]
\addplot[color=blue, line width = 0.5mm, dashed]
coordinates{(0.1, 0.77617377139617*0.239006)(0.01, 0.11110712141108348*0.239006)(0.001,0.011083702841702348*0.239006)(0.0001, 0.0011183542825827786*0.239006)(0.00001, 0.0002321068251407432*0.239006)};
\addlegendentry{\Large{Max Diff}}
\addplot[color=blue,  line width = 0.5mm]
coordinates{(0.1, 0.3783229120175815*0.239006)(0.01, 0.048537280360312324*0.239006)(0.001, 0.00480583667100524*0.239006)(0.0001, 0.00048547774604344187*0.239006)(0.00001, 0.0002552578981875939*0.239006)};
\addlegendentry{\Large{RMSD} }
\end{loglogaxis}
\end{tikzpicture}
\end{subfigure}
\begin{subfigure}[t!]{0.49\textwidth}
 \centering
 \begin{tikzpicture}[scale=0.6]
 \pgfplotsset{
      scale only axis,
  }
\begin{loglogaxis}[
    legend pos=south west, xlabel = \Large{Step size (\AA)},ylabel=\Large{Gradient ($\mathrm{kcal\ mol}^{-1}\mbox{\AA}^{-1}$)},
    legend cell align ={left},
    x dir=reverse, title=\Large{1etn},
    legend style={nodes={scale=0.75, transform shape}}]
\addplot[color=blue, line width = 0.5mm, dashed]
coordinates{(0.1, 1.7752884815201573*0.239006)(0.01, 0.2038530764452945*0.239006)(0.001, 0.09434903278409568*0.239006)(0.0001, 0.0016189198461762366*0.239006)(0.00001, 0.00017148158591639984*0.239006)};
\addlegendentry{\Large{Max Diff}}
\addplot[color=blue,  line width = 0.5mm]
coordinates{(0.1, 0.3357936377933055*0.239006)(0.01, 0.044842366943884035*0.239006)(0.001, 0.00906648627018726*0.239006)(0.0001, 0.00033573815181759276*0.239006)(0.00001, 3.8375485385430035e-05*0.239006)};
\addlegendentry{\Large{RMSD} }
\end{loglogaxis}

\end{tikzpicture}
\end{subfigure}
\caption{Comparison between numerical forces at various step sizes, and the analytical forces for the 1ay3 (left) and 1etn (right) molecules. The two curves report the maximum difference and the root-mean-squared deviation (RMSD) between the two sets of forces.}
\label{fig:numerical_forces}
\end{figure}

Fig.~\ref{fig:numerical_forces} illustrates the convergence of 
the maximum ($\ell^\infty$-error) and the root-mean-squared deviation (RMSD), or equivalently the $\ell^2$-error, of the error vector $\mathsf{Err}$ as a function of $h$ and first-order convergence is indeed observed.
However, beyond $h=10^{-5}$, the finite precision of the algorithms interferes with the convergence of the numerical forces.
We deduce correctness of our theory and implementation from these tests.




\subsubsection{Accuracy of Energy and Forces}

As a preliminary test, we now investigated the role of the discretization parameters used in ddLPB. We present results for the four smallest structures (1ay3, 1etn, 1du9, 1d3w). For 1ay3, 1etn, and 1du9 we run the calculations without the FMM acceleration, for 1d3w we use the FMM acceleration using $p_{\max}$ = 20.  We run a series of energy and force calculations using different values of $\ell_{\max} ( =  2, 4, 6, 8, 10, 12)$.

These calculations were run using a tight convergence threshold of $10^{-8}$ for the linear system, and a value of $N_{\leb}$ = 590, which is enough to perform the numerical quadrature of the high order spherical harmonics used in this test. 

\begin{figure}[ht]
    \centering
    \includegraphics[width=8cm]{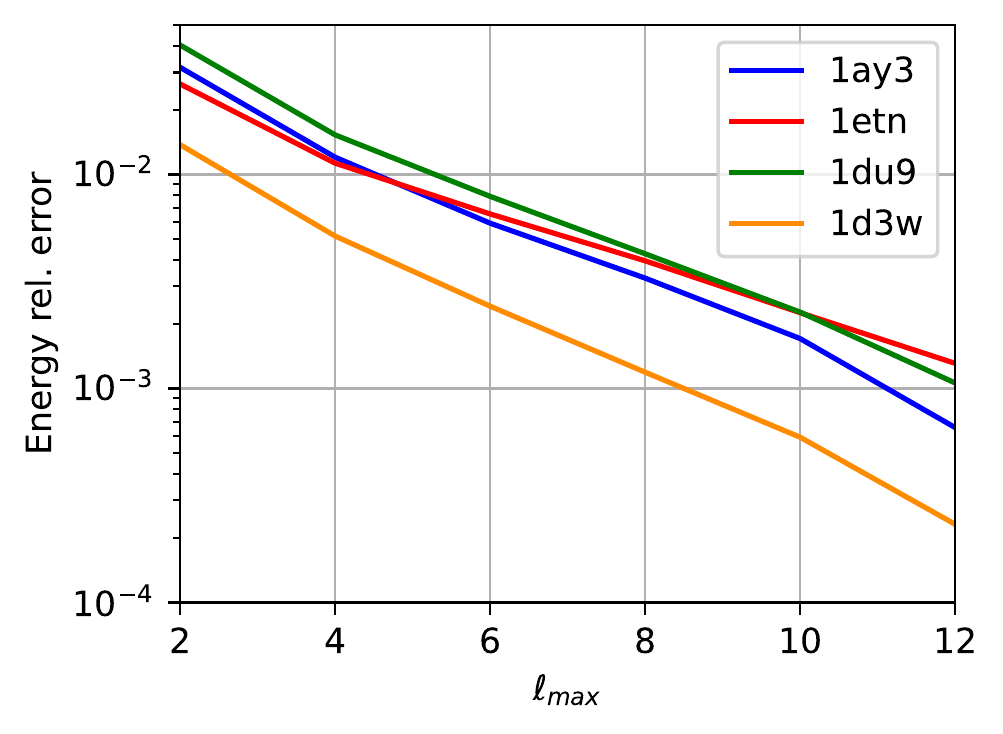}
    \caption{Relative error of the energy against the reference values, with respect to the discretization parameter $\ell_{\max}$. The analysis has been done for the four smallest structures. The reference values are obtained through an exponential fitting of the available energies.}
    \label{fig:lmax_benchmark}
\end{figure}

For each structure, we estimated the converged energy by first using an exponential fitting of the energy with respect to $\ell_{\max}$, and then taking the limit $\ell_{\max} \rightarrow \infty$. Once the reference values were available, we computed the relative error of the discretization at each value of $\ell_{\max}$, for each structure. These results are shown in Fig.~\ref{fig:lmax_benchmark}. 
We draw two conclusions from these tests. 
First, we observe that it is possible to achieve an accuracy below 1\% of the energy by taking $\ell_{\max}$ = 6.
Second, we observe an exponential decay of the error with respect to $\ell_{\max}$ (which justifies the exponential fitting). This allows to reach relatively quickly a regime of high-accuracy with a moderate number of degrees of freedom.

Next, we investigate the role of the parameter $p_{\max}$ which is used to control the FMM acceleration and accuracy.  For each value of $\ell_{\max}$, and for each of the three smallest structures, the reference value is obtained with a non FMM calculation. For the 1d3w structure, the reference value is obtained with an FMM calculation for $p_{\max} = 20$. For each structure, and for each value of $\ell_{\max}$ we run a series of calculations using different values of $p_{\max} ( = 2, 4, 6, 8, 10, 12) $. Also in this case, the convergence threshold was set to $10^{-8}$, and $N_{\leb} = 590$.


\begin{figure}[ht]
    \centering
    \includegraphics[width=14cm]{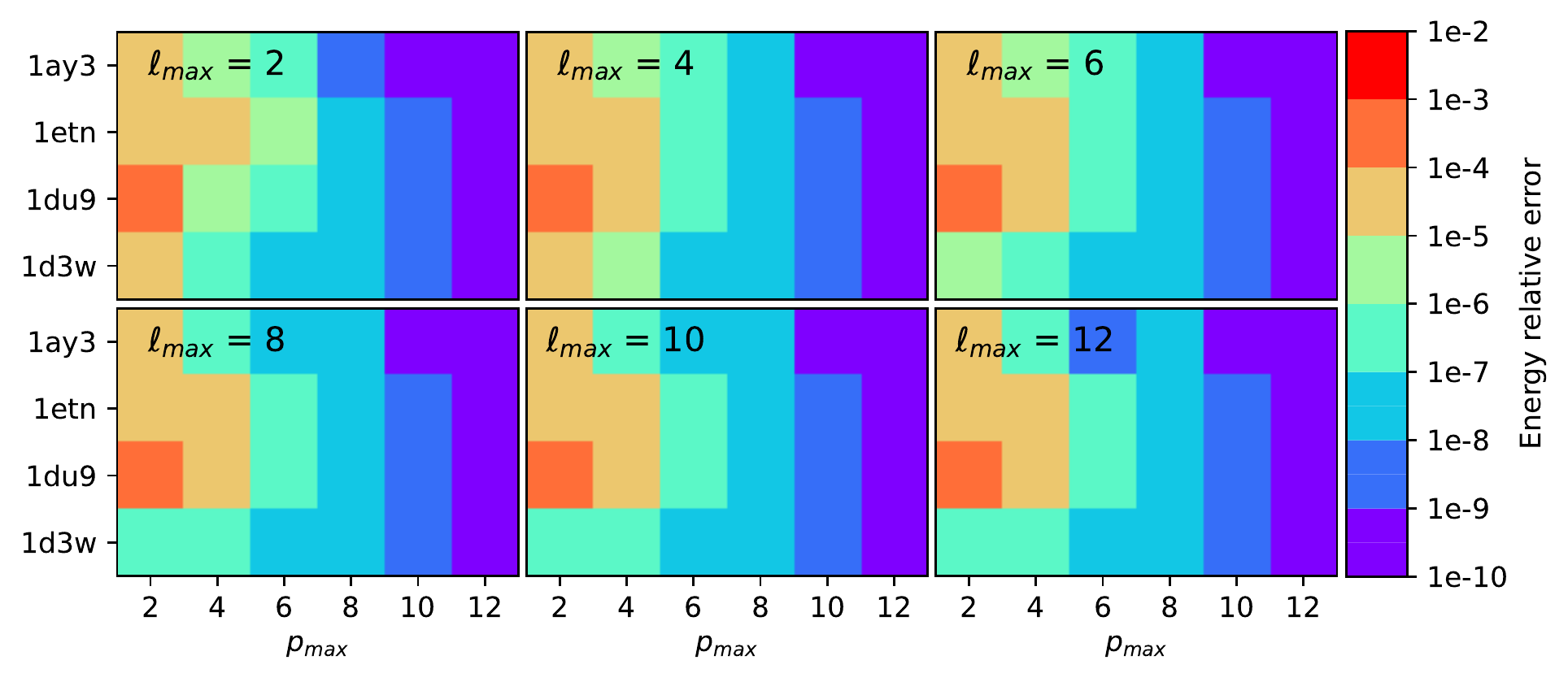}
    \caption{Absolute error of the energy against the reference values, with respect to the FMM discretization parameter $p_{\max}$. The analysis has been done for the four smallest structures, for different values of $\ell_{\max}$ values.}
    \label{fig:fmm_energy_benchmark}
\end{figure}

\begin{figure}[ht]
    \centering
    \includegraphics[width=14cm]{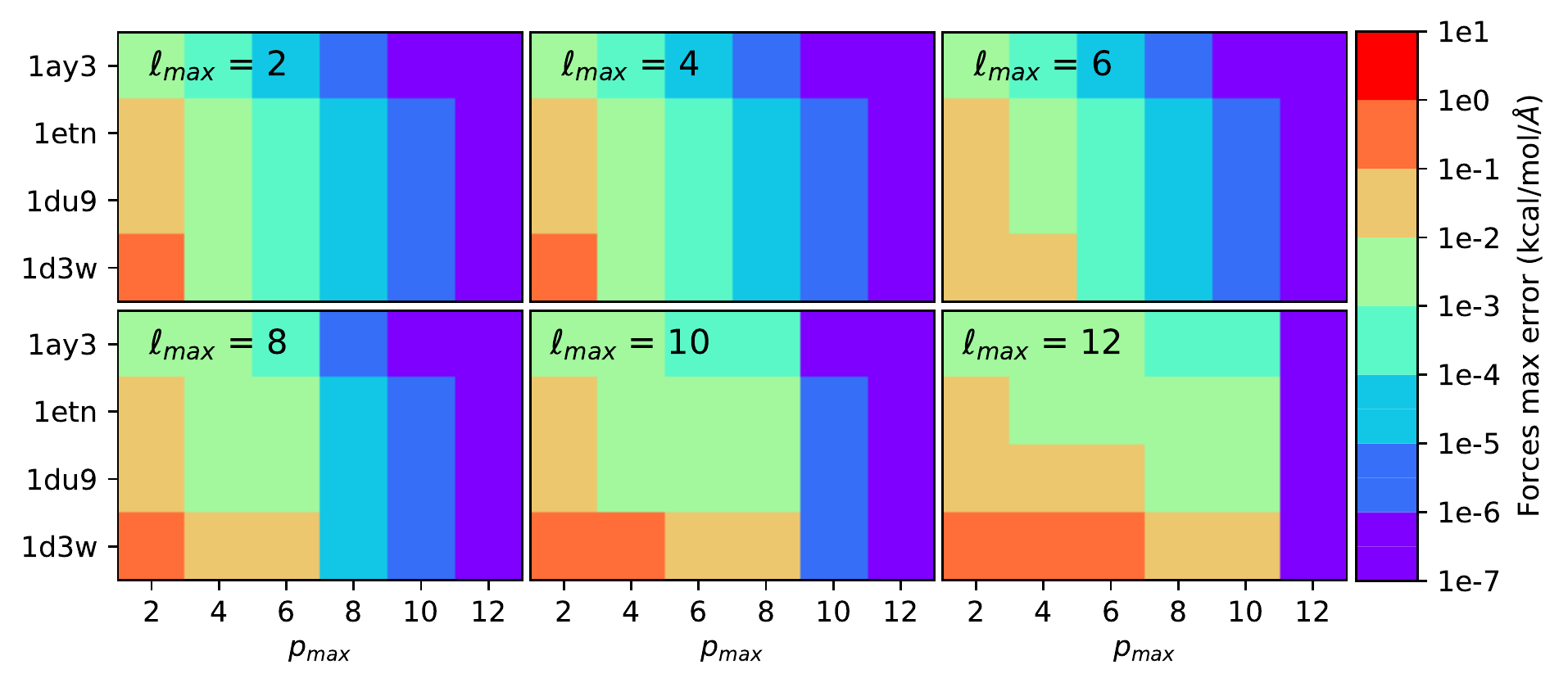}
    \caption{Maximum absolute error of the forces against the reference values, with respect to the FMM discretization parameter $p_{\max}$. The analysis has been done for the four smallest structures, for different values of $\ell_{\max}$ values. }
    \label{fig:fmm_forces_benchmark}
\end{figure}

For each calculation we computed the relative error on the energy using as a reference the corresponding non FMM accelerated calculation (or in case of the system 1d3w, the FMM calculations obtained by setting $p_{\max} = $ 20). Furthermore, we also computed the maximum error on the forces. These results are plotted in Figs.~\ref{fig:fmm_energy_benchmark}~and~\ref{fig:fmm_forces_benchmark}. The plots show that the energy is particularly robust with respect to the FMM discretization, however, the forces need a value of $p_{\max}$ at least equal to $\ell_{\max}$ to achieve a high accuracy,  the same observation was found in the recent publication \cite{mikhalev_nottoli_stamm_2022}.

\subsubsection{Rotational Symmetry of the ddLPB Model}
\begin{figure}[t!]
    \centering
\begin{tikzpicture}[scale=0.5]
\begin{axis}[
    legend pos=north east, xlabel = \Large{$\theta$}, ylabel= \Large{Energy ($\mathrm{kcal\ mol}^{-1}$)}, ymin = -2.19*0.239006, ymax = -2.11*0.239006,
    ytick={-0.525,-0.52,-0.515, -0.51, -0.505, -0.5},
    yticklabels={-0.525,-0.52,-0.515, -0.51, -0.505, -0.5},
    legend cell align ={left}, title = {\Large{$\ell_{\max} = 6$, $N_{\leb} = 110$}}]
\addplot[color=blue,  mark=square*, line width = 0.5mm, dashdotted,, mark options = {scale= 1.0, solid}] 
coordinates{(  0.000000
 , -2.141167484175086*0.239006 )(  0.153248
 , -2.128158655977242*0.239006 )(  0.306497
 , -2.1219648937279962*0.239006 )(  0.459745
 , -2.1363631700005303*0.239006 )(  0.612994
 , -2.1233180994430545*0.239006 )(  0.766242
 , -2.1869849105705383*0.239006 )(  0.919491
 , -2.138394293868854*0.239006 )(  1.072739
 , -2.127490116519608*0.239006 )(  1.225987
 , -2.126476601331254*0.239006 )(  1.379236
 , -2.130953387730211*0.239006 )(  1.532484
 , -2.1390630004889695*0.239006 )(  1.685733
 , -2.126905577429699*0.239006 )(  1.838981
 , -2.1222781718091146*0.239006 )(  1.992229
 , -2.137297269760497*0.239006 )(  2.145478
 , -2.112966067995145*0.239006 )(  2.298726
 , -2.171710147250639*0.239006 )(  2.451975
 , -2.1541073529470736*0.239006 )(  2.605223
 , -2.115744641000984*0.239006 )(  2.758472
 , -2.1330235156201014*0.239006 )(  2.911720
 , -2.1275425981071225*0.239006 )(  3.064968
 , -2.132637891149738*0.239006 )(  3.218217
 , -2.1326378911487596*0.239006 )(  3.371465 
 , -2.127542598131357*0.239006 )(  3.524714 
 , -2.133023515597977*0.239006 )(  3.677962 
 , -2.1157446410010543*0.239006 )(  3.831211 
 , -2.1541073529470744*0.239006 )(  3.984459 
 , -2.1717101472506464*0.239006 )(  4.137707 
 , -2.1129660679943956*0.239006 )(  4.290956 
 , -2.137297269760496*0.239006 )(  4.444204 
 , -2.1222781718192483*0.239006 )(  4.597453 
 , -2.126905577428904*0.239006 )(  4.750701 
 , -2.1390630004772473*0.239006 )(  4.903950 
 , -2.130953387730219*0.239006 )(  5.057198 
 , -2.1264766013311527*0.239006 )(  5.210446 
 , -2.1274901165195925*0.239006 )(  5.363695 
 , -2.1383942938435827*0.239006 )(  5.516943 
 , -2.186984910467808*0.239006 )(  5.670192 
 , -2.123318099443055*0.239006 )(  5.823440 
 , -2.1363631700005334*0.239006 )(  5.976688
 , -2.1219648936958406*0.239006 )(  6.129937
 , -2.128158655977242*0.239006 )};
\end{axis}
\end{tikzpicture}
\begin{tikzpicture}[scale=0.5]
\begin{axis}[
    legend pos=north east, xlabel = \Large{$\theta$}, ylabel= \Large{Energy ($\mathrm{kcal\ mol}^{-1}$)}, ymin = -2.19*0.239006, ymax = -2.11*0.239006,
    ytick={-0.525,-0.52,-0.515, -0.51, -0.505, -0.5},
    yticklabels={-0.525,-0.52,-0.515, -0.51, -0.505, -0.5},
    legend cell align ={left}, title = {\Large{$\ell_{\max} = 10$, $N_{\leb} = 302$}}]
\addplot[color=blue,  mark=square*, line width = 0.5mm, dashdotted,, mark options = {scale= 1.0, solid}] 
coordinates{(  0.000000
 , -2.142153506849608*0.239006 )(  0.153248
 , -2.140654993228323*0.239006 )(  0.306497
 , -2.1389837500345648*0.239006 )(  0.459745
 , -2.133516686089939*0.239006 )(  0.612994
 , -2.126706390107559*0.239006 )(  0.766242
 , -2.1538272500686664*0.239006 )(  0.919491
 , -2.125722442374043*0.239006 )(  1.072739
 , -2.1358723121720513*0.239006 )(  1.225987
 , -2.140646127998721*0.239006 )(  1.379236
 , -2.1408486230121175*0.239006 )(  1.532484
 , -2.139851716984284*0.239006 )(  1.685733
 , -2.1382872252065597*0.239006 )(  1.838981
 , -2.1382927962873106*0.239006 )(  1.992229
 , -2.134525709882732*0.239006 )(  2.145478
 , -2.1328939293597333*0.239006 )(  2.298726
 , -2.1461647060087956*0.239006 )(  2.451975
 , -2.134031937120269*0.239006 )(  2.605223
 , -2.1372031532169675*0.239006 )(  2.758472
 , -2.1389230316655787*0.239006 )(  2.911720
 , -2.139478398863097*0.239006 )(  3.064968
 , -2.136958693755489*0.239006 )(  3.218217
 , -2.1369586937583724*0.239006 )(  3.371465 
 , -2.1394783988630985*0.239006 )(  3.524714 
 , -2.138923031665579*0.239006 )(  3.677962 
 , -2.1372031532169666*0.239006 )(  3.831211 
 , -2.1340319371605316*0.239006 )(  3.984459 
 , -2.1461647060087965*0.239006 )(  4.137707 
 , -2.132893929380476*0.239006 )(  4.290956 
 , -2.134525709882732*0.239006 )(  4.444204 
 , -2.138292796285142*0.239006 )(  4.597453 
 , -2.1382872252065597*0.239006 )(  4.750701 
 , -2.1398517169842823*0.239006 )(  4.903950 
 , -2.1408486230121153*0.239006 )(  5.057198 
 , -2.14064612799872*0.239006 )(  5.210446 
 , -2.1358723122090515*0.239006 )(  5.363695 
 , -2.1257224423740437*0.239006 )(  5.516943 
 , -2.15382725006867*0.239006 )(  5.670192 
 , -2.12670639010756*0.239006 )(  5.823440 
 , -2.1335166860899513*0.239006 )(  5.976688
 , -2.138983750034566*0.239006 )(  6.129937
 , -2.140654993228319*0.239006 )};
\end{axis}
\end{tikzpicture}
\begin{tikzpicture}[scale=0.5]
\begin{axis}[
    legend pos=north east, xlabel = \Large{$\theta$}, ylabel= \Large{Energy ($\mathrm{kcal\ mol}^{-1}$)}, ymin = -2.19*0.239006, ymax = -2.11*0.239006,
    ytick={-0.525,-0.52,-0.515, -0.51, -0.505, -0.5},
    yticklabels={-0.525,-0.52,-0.515, -0.51, -0.505, -0.5},
    legend cell align ={left}, title = {\Large{$\ell_{\max} = 15$, $N_{\leb} = 770$}}]
\addplot[color=blue,  mark=square*, line width = 0.5mm, dashdotted,, mark options = {scale= 1.0, solid}] 
coordinates{(  0.000000
 , -2.136527080778879*0.239006 )(  0.153248
 , -2.1367460746122284*0.239006 )(  0.306497
 , -2.1373619760131852*0.239006 )(  0.459745
 , -2.135909541056128*0.239006 )(  0.612994
 , -2.1358447844876336*0.239006 )(  0.766242
 , -2.1356677507388957*0.239006 )(  0.919491
 , -2.1383600542569554*0.239006 )(  1.072739
 , -2.1361238250966195*0.239006 )(  1.225987
 , -2.1370157634435083*0.239006 )(  1.379236
 , -2.136919046474587*0.239006 )(  1.532484
 , -2.1365345997120784*0.239006 )(  1.685733
 , -2.136475447644429*0.239006 )(  1.838981
 , -2.1371106954368764*0.239006 )(  1.992229
 , -2.1370662858863985*0.239006 )(  2.145478
 , -2.135308643026792*0.239006 )(  2.298726
 , -2.1362000016216074*0.239006 )(  2.451975
 , -2.138042996669981*0.239006 )(  2.605223
 , -2.13665927152422*0.239006 )(  2.758472
 , -2.1370127359059854*0.239006 )(  2.911720
 , -2.136683564748673*0.239006 )(  3.064968
 , -2.136673614919522*0.239006 )(  3.218217
 , -2.1366736149195096*0.239006 )(  3.371465 
 , -2.1366835647489872*0.239006 )(  3.524714 
 , -2.137012735905985*0.239006 )(  3.677962 
 , -2.136659271525697*0.239006 )(  3.831211 
 , -2.138042996656149*0.239006 )(  3.984459 
 , -2.136200001621645*0.239006 )(  4.137707 
 , -2.1353086430267947*0.239006 )(  4.290956 
 , -2.137066285877199*0.239006 )(  4.444204 
 , -2.13711069543688*0.239006 )(  4.597453 
 , -2.136475447665295*0.239006 )(  4.750701 
 , -2.1365345997182286*0.239006 )(  4.903950 
 , -2.136919046473843*0.239006 )(  5.057198 
 , -2.137015763443499*0.239006 )(  5.210446 
 , -2.1361238250965133*0.239006 )(  5.363695 
 , -2.138360054256841*0.239006 )(  5.516943 
 , -2.135667750738892*0.239006 )(  5.670192 
 , -2.135844784487631*0.239006 )(  5.823440 
 , -2.1359095410545583*0.239006 )(  5.976688
 , -2.1373619760131546*0.239006 )(  6.129937
 , -2.1367460746094507*0.239006 )};
\end{axis}
\end{tikzpicture}
\caption{Rotational symmetry for Hydrogen Fluoride molecule, for $\ell_{\max} =6, N_{\leb} = 110$; $\ell_{\max} = 10, N_{\leb}=302,$ and $\ell_{\max}=15, N_{\leb}=770$. $\theta$ corresponds to the angle between the atoms in a fixed reference system.}
\label{fig:rotational_symmetry}
\end{figure}
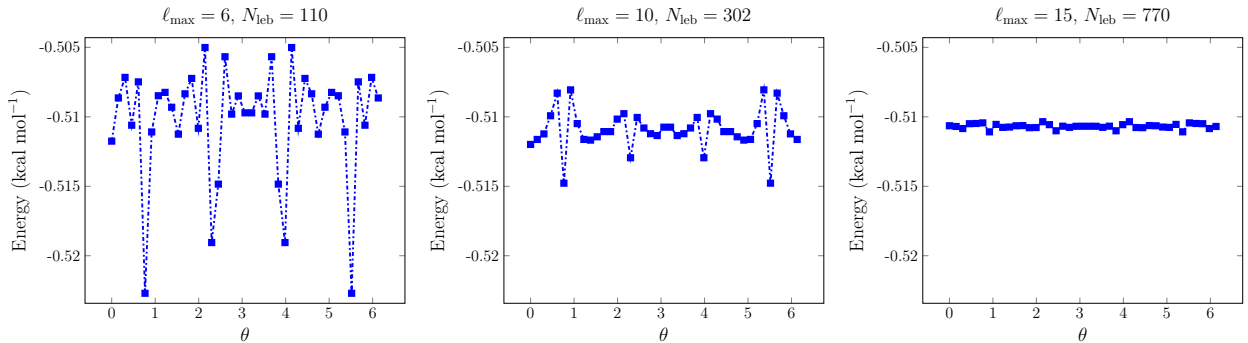

In this example we show that the fluctuation of the energy computation of the ddLPB model  under rotational symmetry is systematically controllable.
For this we use the Hydrogen Fluoride molecule and run the simulations with keeping the Hydrogen atom fixed at $(0,0,0)$ and rotating the Fluorine atom around the Hydrogen atom with $\theta\in[0,2\pi]$, where $\theta$ is the angle subtended by the center of Fluorine and Hydrogen atom. To obtain accurate quadrature we set the number of quadrature points propotional to the degree of spherical harmonics as given in \cite{CMS13}. Fig.~\ref{fig:rotational_symmetry} presents the energy for values of $\ell_{\max}=6,10,15$. We notice that the energy fluctuation under rotation of the fluorine atom is systematically controllable by the discretization parameter $\ell_{\max}$ and is about 0.017 $\mathrm{kcal\ mol}^{-1}$ for the coarsest discretization that is presented here.

\subsection{Complexity of the Discretisation}

\subsubsection{Scaling of ddLPB}

After finding suitable parameters for achieving the required accuracy for the ddLPB-energy, in this section we investigate the performance of the method. To do this, an energy and force calculations were run for each structure. For these calculations the parameters were set as $\ell_{\max} = 6$, $N_{\leb} = 110$, $\mathtt{tol}= 10^{-4}$, and $p_{\max} = 6$. We run these calculations with two setups, in one case by computing the sparse matrix vector products ($AX_r$ and $B X_e$) ``onthefly'' (i.e. without assembling the matrices), and, in the other case, by storing the sparse matrices ``incore'' and using \textsc{BLAS} routines to perform the matrix vector products. In principle the second strategy should be faster but at the cost of an increased, but still linear scaling, memory usage.


\begin{figure}[ht]
    \centering
    \includegraphics[width=8cm]{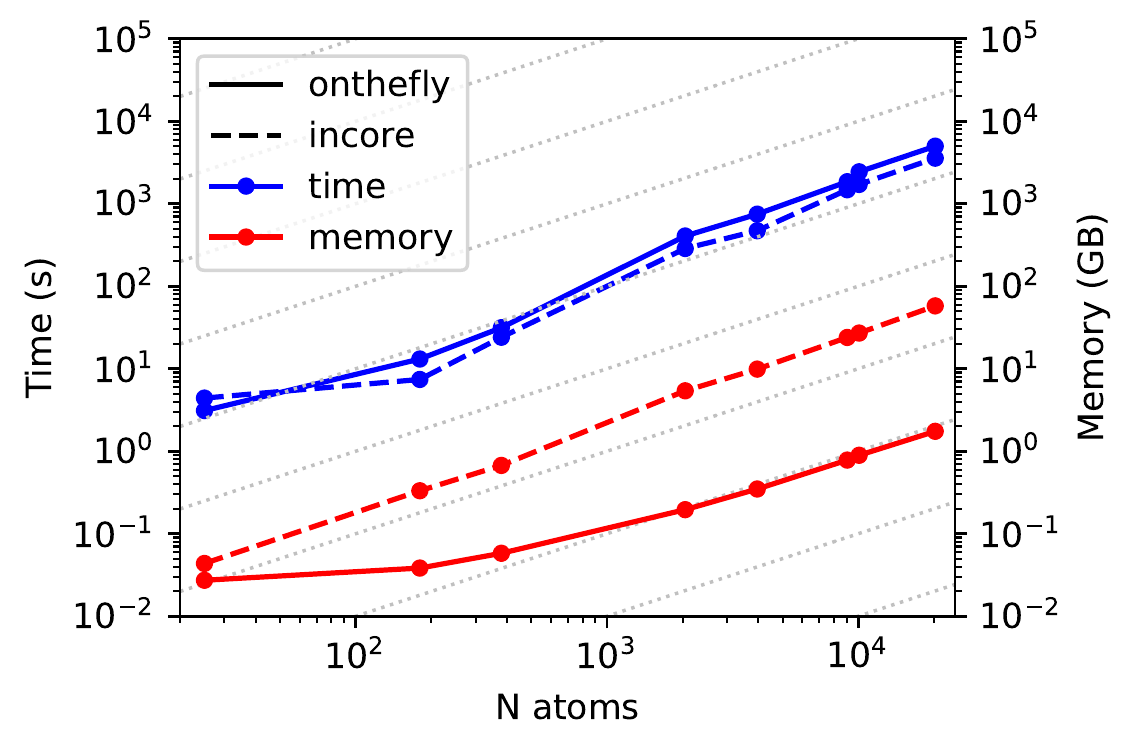}
    \caption{Time and memory required to perform a complete ddLPB force calculation for different structures with $\ell_{\max}=6,N_{\leb}=110,\mathtt{tol}=10^{-4},$ and $p_{\max}=6$. The slope corresponding to a linear scaling regime is highlighted with gray dotted lines. The ``incore'' results are reported as dashed lines, the ``onthefly'' results as solid lines.}
    \label{fig:time_memory}
\end{figure}

Time and memory required by each calculation were measured using the Unix program \texttt{time} and the results are plotted in Fig.~\ref{fig:time_memory}. The log--log plot confirms that the ddLPB method is linear scaling in both time and memory. The linear scaling regime is retained in both the ``incore'' and ``onthefly'' setup, however the ``incore'' setup is only slightly faster than the ``onfly'' setup but at the cost of a significantly increased memory usage.

\begin{figure}[ht]
    \centering
    \includegraphics[width=8cm]{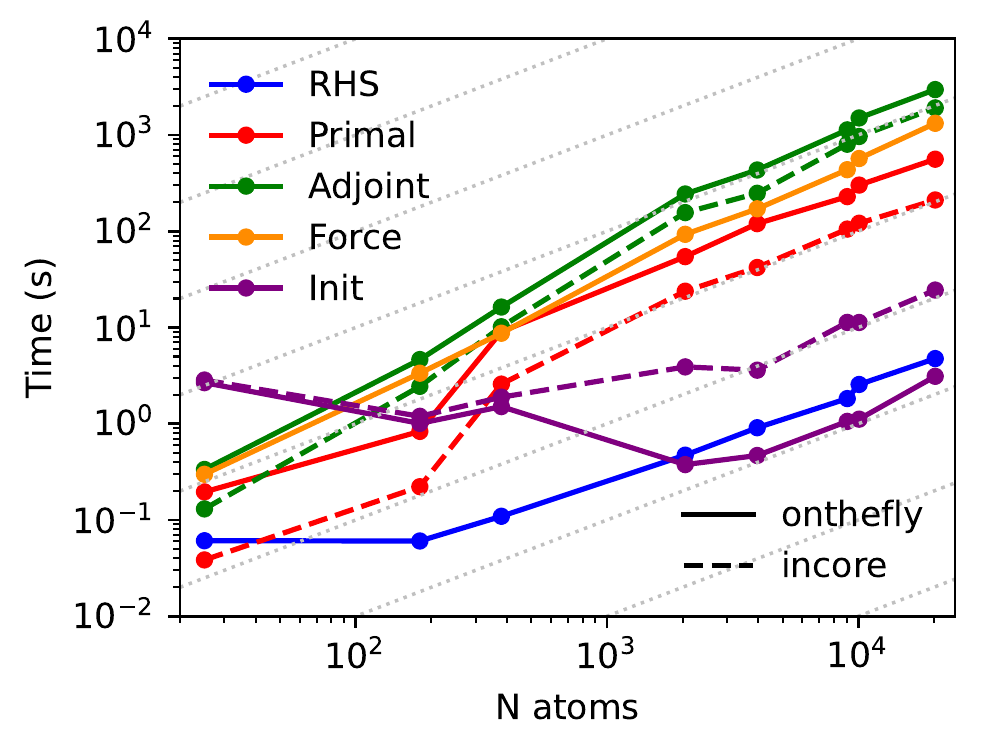}
    \caption{Time and memory required to perform a complete ddLPB force calculation for different structures with $\ell_{\max}=6,N_{\leb}=110,\mathtt{tol}=10^{-4},$ and $p_{\max}=6$. The slope corresponding to a linear scaling regime is highlighted with gray dotted lines. For the linear systems, the ``incore'' results are reported as dashed lines, the ``onthefly'' results as solid lines.}
    \label{fig:steps_time}
\end{figure}

In Fig.~\ref{fig:steps_time} we report a detailed breakdown of the time required to perform the various steps of the ddLPB calculation.
The initialization time, as well as the time required to compute the RHS (electric potential, electric field, and electric field gradient) are linear scaling and negligible with respect to the rest. The time required to solve the primal and the adjoint linear systems are two of the main contributions, both of them are linear scaling, and are slightly faster when the ``incore'' setup is used. Finally, the computation of the forces is again linear scaling in time and is of the same order of magnitude as of solving the linear systems.

\subsubsection{Comparison with Other Software}\label{sec:comparison_software}

Next, we compare the ddLPB model with some of the well used software, namely APBS-FDM and TABI-PB.

For APBS-FDM, the calculations were performed using the box provided by PDB2PQR (keyword \texttt{key}), which is enough to contain the structures, and a number of grid points (keyword \texttt{grid}) suitable for the multigrid algorithm, calculated using
\begin{equation}
n = c \, 2^{\ell + 1} + 1,
\end{equation}
where $n$ is the number of grid points along a given dimension, $\ell$ is the depth of the multilevel solver (keyword \texttt{nlev}), and $c$ is an arbitrary integer. We choose $c$ such that, with $\ell = -4$, a certain target density of points is achieved. In the following discussion, we report the actual density of points computed with $\ell = -4$ and as an average over the three dimensions.
The remaining relevant keywords are \texttt{chgm = spl4}, \texttt{bcfl = mdh}, \texttt{srad = 0.0}, and \texttt{swin = 0.3}.
For ddLPB, we set the tolerance $\texttt{tol} = 10^{-4}$, $N_{\leb}=590$, and the smooth-switching window $\eta = 0.0$. The maximum degree of spherical harmonics is set to values between 2 and 12. Furthermore, these calculations were run using both the ``incore'' and ``onthefly'' setup.
Finally, for TABI-PB we used the following keywords: \texttt{mesh = SES}, \texttt{pdie = 1}, \texttt{sdie = 78.54}, \texttt{bulk = 0.1}, \texttt{temp = 298.15}, \texttt{tree\_degree = 2}, \texttt{tree\_max\_per\_leaf = 50}, and \texttt{tree\_theta = 0.8}. The density of points (\texttt{sdens}) was varied from 5 to 40 to study the convergence of the results. All the calculations were run using a single core, to ease further comparisons.

Once a series of points for different discretizations were gathered, we performed a fitting to extrapolate the energy value in the limit of an infinite discretization. It can be shown that the energies computed using both TABI-PB and APBS-FDM converge in an algebraic way with respect to the number of degrees of freedom, on the other hand, the energies computed using ddLPB converge exponentially with respect to the number of degrees of freedom. For this reason, for the first method we use a nonlinear fitting of the form $a + bx^c$ where $x$ is the number of triangles of the cavity; for the second model we use a linear fitting $a + b x$ where $x$ is the inverse of the average grid spacing in~{\AA}; and finally for ddLPB we use an exponential fitting of the form $a + b e^{c x}$ where $x$ is the maximum degree of the spherical harmonics. In each case, $a$ is the extrapolated energy in the limit of an infinite discretization, which was used to compare the results of ddLPB and APBS-FDM, and to compute the discretization errors.

\begin{table}[ht]
    \centering
    \begin{tabular}{|crrr|}
    \hline
    \textbf{Structure} & \textbf{ddLPB} & \textbf{APBS-FDM} & \textbf{Rel. diff. (\%)} \\ \hline\hline
1ay3 & -31.2 & -31.4 & 0.47 \\
1etn & -126.7 & -125.6 & 0.91 \\
1du9 & -296.8 & -295.6 & 0.39 \\
1d3w & -3384.3 & -3358.5 & 0.77 \\
1jvu & -1563.8 & -1555.8 & 0.52 \\
 \hline
    \end{tabular}
    \caption{Comparison of the ddLPB and APBS-FDM energies (kcal mol$^{-1}$) in the extrapolated limit of an infinite discretization. The last column reports the percent relative differences computed as {$|\text{ddLPB} - \text{APBS-FDM}|/|\text{APBS-FDM}|\times 100$}.}
    \label{tab:apbs_ddlpb}
\end{table}

Table~\ref{tab:apbs_ddlpb} compares the energies obtained from APBS-FDM and ddLPB for the molecules presented in Table~\ref{tab:structures}.

\begin{figure}
    \centering
    \includegraphics[width=\textwidth]{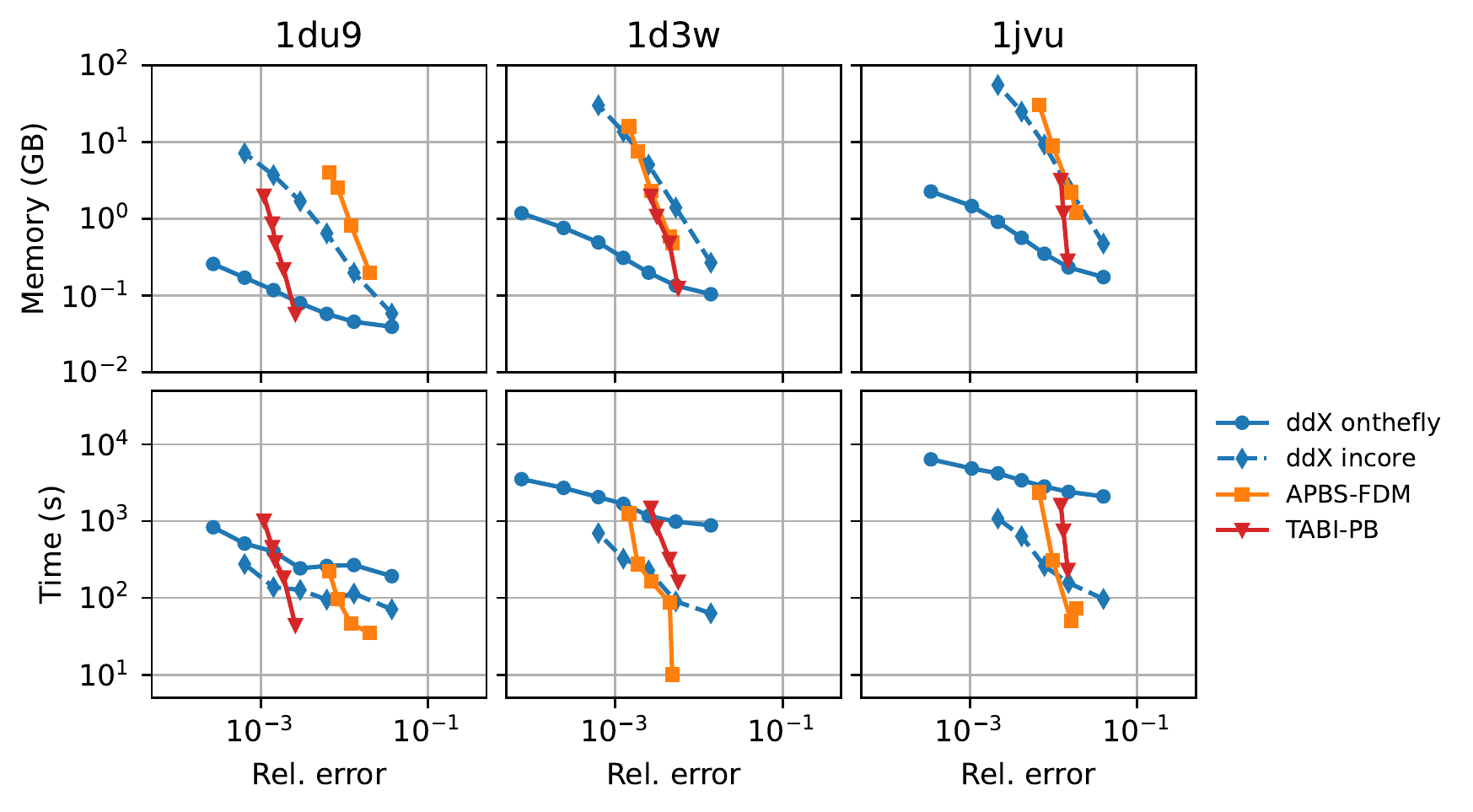}
    \caption{Comparison of the resource consumption for different target accuracies between ddX, APBS-FDM, and TABI-PB.}
    \label{fig:comparison}
\end{figure}

Finally, we present a comparison between the resource consumption of the three methods in Fig.~\ref{fig:comparison}. Since the resource consumption strongly depends on the used discretization, we decided to plot the resources with respect to the discretization error. For this analysis we used only three systems of intermediate size, for which the resource consumption is considerable, but not too large to prevent going to high discretization values.
We observe that, within this computational protocol outlined above and for these molecules, the three methods behave similarly for low accuracy while the exponential convergence makes a real difference if one is aiming for high-accuracy solutions.

\section{Conclusion}\label{sec:summary}
In this work, we provide the detailed derivation of analytical forces and linear scaling for the computation of energy and forces for the ddLPB numerical method which efficiently approximates solutions to the linearized Poisson-Boltzmann equation that is a frequent model used in computational (bio-) chemistry. 
The derivation is technical but mandatory and is based on an adjoint method to compute analytical derivatives of the energy with respect to (possibly many) external parameters such as the nuclear coordinates which result in the computation of the forces.
The implementation of the energy and forces have been validated by a series of benchmark problems and by comparing the results with those of the APBS-FDM-package and TABI-PB.
The current implementation scales linearly with respect to the number of atoms using the fast multipole method (FMM) developed in \cite{mikhalev_nottoli_stamm_2022}.

\appendix
\section{Appendix}\label{appendix_a}
In this appendix we give a brief overview of the fast multipole method (FMM) that was used in Sec.~\ref{sec:numres}. The basic idea of FMM is to reduce the bottleneck, quadratic scaling operations to linear scaling. The quadratic scaling bottlenecks in the computation of energy are the matrix-vector multiplication in the operations corresponding to the primal solution, i.e., Eq.~\eqref{eq:AXg}. For the computation of forces there are two more bottlenecks, namely the matrix-vector multiplication in the computation of the adjoint solution~\eqref{eq:adjoint_system} and the contraction of derivatives described in~\eqref{eq:force_lpb}.
For all these operations, the quadratic scaling is due to the presence of the single layer potential in the nonlocal coupling condition given by~\eqref{eq:nonlocal-coupling} which is affecting only the matrices $\bC_1$, $\bC_2$ and the right hand side $\bF_0$. 

In \cite{mikhalev_nottoli_stamm_2022}, the idea of FMM was introduced for the ddPCM model which uses the Coulomb potential, in contrast to the ddLPB which is based on the Yukawa potential.
We therefore follow the same implementation as proposed in \cite{mikhalev_nottoli_stamm_2022} using a binary adaptive tree-structure, but with adapted multipole-to-multipole (M2M), multipole-to-local (M2L), and local-to-local (L2L) operators for the Yukawa potential.


Indeed, these operations only need to be defined along the OZ-axis and therefore we only report the corresponding OZ translations, i.e., a translation of length $\rho$ along the direction $e_z$.
They are based on notations of equations~(29), (31) and~(32) of the ddPCM-FMM paper and read as

\begin{equation}
    \label{eq:m2m_oz}
    [\mathrm{M2M}(\rho e_z,R_{S}, R_{T})]_{\ell m}^{\ell' m'}
    =
    \begin{cases}
    0, & m \neq m'\\
    \dfrac{ C_{N,\ell m}^{\ell'} \tk_{\ell'}(R_T) }{ \tk_{\ell} (R_{S}) }
    \sum_{k=|m|}^{\min\lbrace \ell',\ell \rbrace} C_{\ell m k}^{\ell'} \frac{ \ti_{\ell+\ell'-k} (\rho) }{ \rho^k }
    , & \text{otherwise},
    \end{cases}
\end{equation}

\begin{equation}
    \label{eq:m2l_oz}
    [\mathrm{M2L}(\rho e_z,R_S, R_T)]_{\ell m}^{\ell' m'}
    =
    \begin{cases}
    0, & m \neq m'\\
   \dfrac{ C_{N,\ell m}^{\ell'} \ti_{\ell'} (R_{T}) (-1)^{\ell} }{ \tk_{\ell} (R_{S}) }
    \sum_{k=|m|}^{ \min\lbrace \ell',\ell \rbrace} C_{\ell m k}^{\ell'} \frac{ \tk_{\ell+\ell'-k} (\rho) }{ (-\rho)^k }
    , & \text{otherwise},
    \end{cases}
\end{equation}

\begin{equation}
    \label{eq:l2l_oz}
    [\mathrm{L2L}(\rho e_z,R_S, R_T)]_{\ell m}^{\ell' m'}
    =
    \begin{cases}
    0, & m \neq m'\\
    \dfrac{ C_{N,\ell m}^{\ell'} \ti_{\ell'} (R_{T}) (-1)^{\ell+\ell'} }{ \ti_{\ell} (R_{S}) }
    \sum_{k=|m|}^{\min\lbrace \ell', \ell \rbrace} C_{\ell m k}^{\ell'} \frac{ \ti_{\ell+\ell'-k} (\rho) }{ \rho^k }
    , & \text{otherwise},
    \end{cases}
\end{equation}
where 
\begin{equation}
C_{N,\ell m}^{\ell'}=\frac{\tilde{N}_{\ell}^m(2\ell'+1)(\ell'-|m|)!(\ell+|m|)!}{\tilde{N}_{\ell'}^{m}}
=C_{N,\ell' m}^{\ell},
\end{equation}
\begin{equation}
C_{\ell mk}^{\ell'}=\frac{(2k)!}{2^k(k+m)!k!(k-m)!(\ell'-k)!(\ell-k)!}=C_{\ell' mk}^{\ell},
\end{equation}
and with 
$\tilde{N}_{\ell}^m$ denoting the normalization factors of spherical harmonics:
\begin{equation}
    \tilde{N}_{\ell}^m =
    \begin{cases}
    \sqrt{ \dfrac{2\ell+1}{4\pi} }, & m = 0\\
    (-1)^m \sqrt{ 2 \cdot \dfrac{2\ell+1}{4 \pi}  \cdot \dfrac{(\ell-|m|)!}{(\ell+|m|)!}}, & m \neq 0.
    \end{cases}
\end{equation}
Note that for the computation of forces, more precisely when the gradients of the potentials and the so-called adjoint potentials (see \cite{mikhalev_nottoli_stamm_2022}) are required, we follow the approach proposed in \cite{mikhalev_nottoli_stamm_2022} which relies on the gradients of the M2M and L2L translations
with identical source and target spheres.
This leads to differentiating OZ-translations of the M2M and L2L operations (equations~\eqref{eq:m2m_oz} and~\eqref{eq:l2l_oz}) with respect to $\rho$ evaluated at $\rho=0$. 
Taking the well-known asymptotic behaviour of $\ti_n(\rho) \approx \frac{\rho^n}{(2n+1)!!}$ near $\rho=0$ and $\ti_0'(0)=0$ into assumption, all the calculations are finally reduced to the following derivative:
\begin{equation}
    \left[ \frac{ \ti_{\ell+\ell'-k} (\rho) }{ \rho^k } \right]'_{\rho=0} 
    =
    \begin{cases}
    0, & \ell+\ell'-2k \neq 1 \\
    \frac{1}{(2k+3)!!}, & \ell+\ell'-2k = 1.
    \end{cases}
\end{equation}
Due to the upper limit $k \leq \min\lbrace \ell, \ell' \rbrace$ the condition $\ell+\ell'-2k=1$ is satisfied only in the case $\ell'=\ell \pm 1$ and $k = \min\lbrace \ell, \ell' \rbrace$.

\section*{Acknowledgements}
CQ is supported by NSFC Grant 12271241, the fund of the Guangdong Provincial Key Laboratory of Computational Science and Material Design (No. 2019B030301001), and Shenzhen Science and Technology Program (No. RCYX20210609104358076).
AJ, AM, and BS are acknowledging support by the German Research Foundation (DFG) under project 440641818.
The authors acknowledge support by the state of Baden-Württemberg through bwHPC.

\bibliographystyle{alpha}
\bibliography{ddLPB}
\end{document}